\documentclass[11pt]{article}
\usepackage{color}
\usepackage[toc,page]{appendix}
\usepackage[utf8]{inputenc}
\usepackage[margin=1in,left=1in]{geometry}
\usepackage[all,knot]{xy}
\usepackage{setspace}
\usepackage{graphicx}
\usepackage{amssymb, amsmath, amsthm, graphics}
\usepackage{mathabx,epsfig}
\usepackage[mathscr]{euscript}
\usepackage{hyperref}
\newcommand{\footremember}[2]{%
    \footnote{#2}
    \newcounter{#1}
    \setcounter{#1}{\value{footnote}}%
}

\usepackage{array, booktabs}  
\usepackage{caption}    

\usepackage{tikz}
\usetikzlibrary{calc}
\usetikzlibrary{matrix}
\usepackage{latexsym}
\usepackage[all]{xy}
\usepackage{color}

\usepackage{enumitem} 

\usepackage{bbold}

\usepackage{tikz}
\input xy
\xyoption{all}
\usepackage{euscript}
\usepackage{multirow}
\usepackage{etex, pictexwd,dcpic}
\usetikzlibrary{positioning}
\usetikzlibrary{shapes.geometric}
\usetikzlibrary{shapes.misc}
\usetikzlibrary{calc}
\usetikzlibrary{positioning}

\usepackage{pgflibraryarrows}
\usepackage{pgflibrarysnakes}

\usetikzlibrary{trees} 
\usetikzlibrary[trees] 

\newtheorem{theorem}{Theorem}[section]
\newtheorem{definition}[theorem]{Definition}
\newtheorem{lemma}[theorem]{Lemma}
\newtheorem{corollary}[theorem]{Corollary}
\newtheorem{proposition}[theorem]{Proposition}

\newtheorem{remark}[theorem]{Remark}

\newtheorem{example}[theorem]{Example}

\numberwithin{equation}{section}

\newcommand{\introthmname}{}
\newtheorem{introthminn}{\introthmname}
\newenvironment{introthm}[1]
  {\renewcommand{\introthmname}{#1}\begin{introthminn}}
  {\end{introthminn}}

\DeclareMathAlphabet\mathbfcal{OMS}{cmsy}{b}{n}

\renewcommand{\(}{\begin{equation*}}
\renewcommand{\)}{\end{equation*}}
\newcommand{\bea}{\begin{eqnarray*}}
\newcommand{\eea}{\end{eqnarray*}}

\makeatletter
\newcommand{\colim@}[2]{%
  \vtop{\m@th\ialign{##\cr
    \hfil$#1\operator@font lim$\hfil\cr
    \noalign{\nointerlineskip\kern1.5\ex@}#2\cr
    \noalign{\nointerlineskip\kern-\ex@}\cr}}%
}
\newcommand{\colim}{%
  \mathop{\mathpalette\colim@{\rightarrowfill@\textstyle}}\nmlimits@
}
\makeatother

\def\endofproof {\hfill{$\Box$}\\}

\def\H{\ensuremath{\ES{H}}}

\def\H{\ensuremath{\ES{H}}}

\newcommand{\beq}{\begin{equation}}
\newcommand{\eeq}{\end{equation}}

\newcommand{\into}{\hookrightarrow}

\newcommand{\ES}[1]{\ensuremath{\EuScript{#1}}}

\newcommand{\theproof}{\noindent {\bf Proof.\ }}

\numberwithin{equation}{section}

\renewcommand{\(}{\begin{equation}}
\renewcommand{\)}{\end{equation}}

\def\1{{\bf 1}}

\def\<{\langle}
\def\>{\rangle}

\numberwithin{equation}{section}


\newcommand{\R}{\ensuremath{\mathbb R}}
\newcommand{\RR}{\ensuremath{\mathbb R}}

\newcommand{\ZZ}{\ensuremath{\mathbb Z}}
\newcommand{\Z}{\ensuremath{\mathbb Z}}

\newcommand{\BB}{\ensuremath{\mathbf B}}
\newcommand{\MM}{\ensuremath{\mathbf M}}
\newcommand{\TT}{\ensuremath{\mathbf T}}
\newcommand{\CC}{\ensuremath{\mathbb C}}

\newcommand{\sset}{\ensuremath{s\mathscr{S}\mathrm{et}}}

\newcommand{\set}{\ensuremath{\mathscr{S}\mathrm{et}}}

\newcommand{\sh}{\ensuremath{\mathscr{S}\mathrm{h}}}

\newcommand{\PSh}{\mathscr{PS}\mathrm{h}_\infty}
\newcommand{\Sh}{\mathscr{S}\mathrm{h}_\infty}
\newcommand{\cartsp}{\mathscr{C}\mathrm{art}\mathscr{S}\mathrm{p}}

\newcommand{\map}{\mathrm{Map}}

\begin{document}

\title{Parametrized geometric cobordism and smooth Thom stacks}

\author{
  Daniel Grady \footremember{alley}{djg9@nyu.edu}
  \and Hisham Sati \footremember{trailer}{hsati@nyu.edu}
  }


\maketitle

\begin{abstract}
We develop a theory of parametrized geometric cobordism by introducing
smooth Thom stacks. This requires identifying and constructing a smooth 
representative of the Thom functor acting on vector bundles equipped with 
extra geometric data, leading to a geometric refinement of the the Pontrjagin-Thom
construction in stacks. We demonstrate that the resulting theory generalizes  the 
parametrized cobordism of Galatius-Madsen-Tillman-Weiss. 
The theory has the feature of being both 
versatile and general, allowing for the inclusion of families of various geometric 
data, such as metrics on manifolds and connections on vector bundles, as
in recent work of Cohen-Galatius-Kitchloo and Ayala. 
\end{abstract} 

\tableofcontents

\section{Introduction and overview}

Bordism was introduced by Pontrjagin in order to understand manifolds
from a homotopy-theoretic point of view (see \cite{Po}). Thom showed that 
cobordism groups  could be computed by means of homotopy theory, via the 
Thom complex construction \cite{Th}. 
Excellent surveys can be found in \cite{Ad}\cite{Pe}\cite{Ru}\cite{St}.
Being originally a geometric construction (see also \cite{Kos} for more on this approach),
it is then natural to aim to tie back homotopic constructions with differential geometric ones.

\medskip
Conner and Floyd \cite[Section I.9]{CF}  define differentiable bordism groups 
$D_n(X^k)$ of degree $n$  of a $k$-dimensional differentiable manifold 
without boundary as follows. 
Consider pairs $(M^n, f)$ consisting of a closed unoriented manifold $M^n$
and a differentiable map $f: M^n \to X^k$. Such a pair 
bords if and only if there is a compact 
manifold $B^{n+1}$ with $\partial B^{n+1}=M^n$ and a differentiable map
$g: B^{n+1} \to X^k$, with 
$g|_{\partial B^{n+1}}=f$. In order for the bordism relation to be transitive, 
it is required further that there 
exist an open set $U\supset \partial B^{n+1}$ and a diffeomorphism 
 $h: \partial B^{n+1} \times [0, 1)\to U$ with $h(x, 0)=x$ and $g(h(x, t))=f(x)$ for 
 all $0 \leq t < 1$ and $x \in \partial B^{n+1}$. 
 The resulting group $D_n(X^k)$ admits a natural homomorphism to the 
 unoriented cobordism group $D_n(X^k) \to \mathfrak{N}_n(X^k)$, 
 taking the class $[(M^n, f)]$ to the class $[(M^n, f)]$, and which is 
 an isomorphism when $X^k$ is a differentiable manifold. 
 The upshot of this is that one can always take the map 
 to be differentiable but then this by itself does not give anything new.

\medskip
Cobordism can be viewed as a generalized cohomology theory \cite{At}\cite{Mil1}. 
One formulation which captures geometric and topological aspects of manifolds is 
differential cohomology, within which developing cobordism is therefore desirable. 
There has been some very foundational recent work in this direction. 
 Hopkins and Singer gave a geometric definition of differential bordism in
the context of differential function spectra in \cite[Section 4.9]{HS}.
 In \cite{BSSW} a differential extension $\widehat{\bf MU}$ of 
complex bordism $MU$ was constructed. The geometry is encoded in cycles, 
requires transversality, hence heavy differential topology, as well as harmonic 
differential characters from $\widehat{\H\Z}$.

\medskip
We provide a general approach to differential cobordism via stacks.
We initially wanted to focus on studying differential refinement of 
unoriented cobordism \cite{Th}\cite{Wa}\cite{Li}\cite{Mit}.
This is the cobordism theory of all compact differentiable manifolds
and, as  such, is perhaps the most general  
cobordism theory, as highlighted by Stong \cite[Chapter VI]{St}. 
From a homotopy point of view, being associated to a $BO$-structure 
makes it more tractable for some purposes than having an additional structure. 
However, our formulation via simplicial sheaves turned out to be general 
enough to allow for differential refinements of cobordisms with $BG$-structures 
in the sense of \cite{St}\cite{Kos}. 

%

\medskip
Our initial motivation for considering cobordism was to study differential refinements of 
cohomology operations. We have constructed primary and secondary differential
cohomology operations in \cite{GS2}\cite{GS1}, and have utilized them in 
constructing spectral sequences in differential cohomology in \cite{GS3}.  
Cohomology operations are intimately related to characteristic classes
and to Thom classes of vector bundles. The way this article has evolved 
\footnote{We started studying differential cobordism as a differential cohomology theory, 
which led us to 
smooth motivic Thom spectrum with connection.  
Generalizing these led us to general smooth motivic Thom spectra 
(different ${\cal F}$'s). While working on this we learned about 
 Madsen-Tillman spectra (as also indicated in acknowledgements) 
 and we started working on the stacky cobordism category and then
  we linked our smooth motivic
construction with that stack via this project.}
made 
it clear that it is better to defer the treatment of characteristic classes to 
a forthcoming article and focus here on development of the smooth cobordism theory
 itself.

\medskip
Thom spaces \cite{Th} provide a convenient setting which involves a smooth structure 
and bundles on $M$.  Note that Hopkins and Singer \cite[Section 4.2]{HS} have defined 
Thom complexes geometrically in the context of differential function spectra,  
and Bunke (\cite[Def 4.181]{Bu1})  has defined differential Thom classes via a slightly 
different approach. The Thom homomorphism has also been considered in other generalized
differential cohomology theories, for instance in differential K-theory by Freed and Lott \cite{FL}. 

\medskip
We seek a smooth analog of a Thom space ${\rm Th}(E)$ of a vector bundle $p:E\to B$. When $p$ is a smooth map between manifolds, this is essentially accomplished by the usual Thom space construction, but with the crucial difference that the quotient is taken not in the ambient category of topological spaces, but rather smooth stacks. As a consequence of this, there are some subtleties which one does not encounter in spaces. For example, recall that in the category of topological spaces 
$\mathscr{T}{\rm op}$, the Thom space of a rank $n$ trivial bundle $p:E\to B$ is 
equivalent to the $n$-fold suspension $\Sigma^nB_+$.
It is an important point for us that this is \emph{not} the case in the category of smooth stacks. 
Essentially, this is because we are remembering the manifold structure on both ${\rm Th}(E)$ 
and $X$, while (the $n$-fold) suspension is a homotopy theoretic construction. We can, 
nevertheless, still be able to move forward  by remembering a little more of the geometry. 
The idea is the following. We observe that if $p:E\to X$ is a trivial vector bundle of rank $n$, defining the Thom 
stack as the quotient of the smooth disc bundle modulo the smooth sphere bundle gives 
an identification
$$
{\rm Th}(E)\simeq D^n/\partial D^n\wedge B_+\;.
$$
This indicates that we should be taking a smooth analogue of the suspension in the category 
of smooth stacks: the quotient stack $D^n/\partial D^n$.  We establish this in Section
\ref{Sec Th}.

\medskip
These considerations seem to point toward an alternative to the usual definition of spectrum in the category of smooth stacks, in which the ``smooth circle" $D^1/\partial D^1$ serves as a model for smooth suspension. 
It turns out that for stacks which are $\RR^1$-invariant (i.e. usual homotopy-invariant) theories, nothing new is gained. However, more geometrically refined cohomology theories, e.g. differential cohomology, are often not homotopy invariant (see e.g. \cite{Bu1} \cite{Urs}) and thus are able to distinguish between the simplicial and smooth spheres. In order to maintain consistency with the goal of this paper, we will not explore the full theory of such objects (which we call \emph{$D^1/\partial D^1$-spectra}) in full generality and save these fundamentals for development elsewhere. Here, we will only provide the definition and focus more on the particular example of differential cobordism in Section \ref{Sec mot}. 

\medskip
In their seminal works \cite{MT}\cite{MW}\cite{GMTW}, Galatius, Madsen, Tillman, and Weiss studied 
the moduli space of embeddings  $j:M\into \RR^{\infty}$ and connected this space with the 
 classifying space of the category of cobordisms and the infinite loop space of the 
 Madsen-Tillman spectrum ${\rm MT}(d)$, with $d$ indexing the dimension of of the bordisms.  
 The homotopy groups of the cobordism  category are given a geometric interpretation in \cite{BSV}. 
 In \cite{GRW} Galatius and Randal-Williams investigated 
subcategories with classifying spaces homotopy equivalent to that of 
the GMTW 
category $C_\theta$ of closed smooth $(d-1)$-manifolds 
and smooth $d$-dimensional cobordisms, equipped with a $\theta$-structure, 
proving a result similar to that of \cite{GMTW}. 
 The homotopy type of the cobordism category 
with objects $(d-1)$-dimensional submanifolds of a fixed 
background manifold $M$ is identified in \cite{OW}. 

\medskip 
 Our goals in this paper are the following. 
 \begin{enumerate}
 
 \vspace{-1mm}
 \item First, to develop a smooth refinement of the Madsen-Tillman spectrum in stacks,  which carries with it geometric data coming from the Thom-space construction. \cite{MT}\cite{GMTW}. 

 \vspace{-3mm}
\item We will then discuss a refined classifying space construction for a smooth category which takes place entirely in the context of smooth stacks. We apply this construction to the smooth cobordism category defined in \cite{GMTW}.

 \vspace{-3mm}
\item Finally, we will describe a refinement of the Pontrjagin-Thom construction 
in this setting and prove an equivalence of two smooth stacks which geometrically realizes to the equivalence proved in \cite{GMTW}. Thus, our construction here can be regarded as a generalization of the main theorem in \cite{GMTW}. 

 \vspace{-3mm}
 \item We will also show how our construction can be used to include connections and geometric data, 
 thus making contact with the setting of Cohen, Galatius and Kitchloo \cite{CGK} and Ayala \cite{A}. 
 \end{enumerate}

We will work in the context of smooth $\infty$-bundles, as presented in 
\cite{NSS}\cite{NSS2}. 
While it might seem overly abstract at first sight, this theory has the desirable advantage of being 
both conceptually simple and unifying, as it recovers various disconnected concepts in the 1-categorical 
setting as being manifested by one single concept in the $\infty$-category context \cite{Lur}. 
However, we will aim to keep abstraction to a minimum, just enough to allows us to present our 
constructions and results. Furthermore, the $\infty$ setting is particularly efficient in keeping track of automorphisms, such as diffeomorphisms or gauge transformations, in a systematic way
(see for instance \cite{FSS} for an illustration). Techniques from $\infty$-categories have been 
used recently in \cite{FSV} to describe higher extensions of diffeomorphism groups  as
group stacks of automorphisms of manifolds, equipped with certain topological structures.
Stacks have also been used in the cobordism context earlier. For instance,  
in \cite{EG} the classical construction of Pontrjagin-Thom maps is extended 
to the category of differentiable local quotient stacks and used to 
detect torsion in the homology of the moduli stack of stable curves.

 \medskip
We present constructions and machinery in stacks in Section \ref{Sec recol}. Most of these have been studied extensively, and from various points of view (see, for instance, \cite{Lur}\cite{DI}\cite{Jar}\cite{Urs}). Although these constructions may seem overly abstract to the reader who is not well versed in abstract homotopy theory, we will see that wading through this level of abstraction  has many advantages. 
Indeed, one of the points of this paper is to show that the classical  theorem of \cite{GMTW} is fairly 
systematic to prove by making use of these techniques. In particular, one does not  need the rather 
complicated transversality theorems used there. These differential topology techniques are even more particularly pronounced in other extensions of the  cobordism  category, for instance to manifolds with corners \cite{Ge}.

 \medskip
Note that whenever we speak of a smooth stack we will mean a simplicial presheaf
on the site of Cartesian spaces with smooth functions between them, i.e., 
an object in (what we called) the category of smooth stacks $\Sh(\cartsp)$ 
(see Definition \ref{smooth stacks}). This might 
be a possible source of confusion, as a smooth stack is usually thought of as a sheaf of 1-categories 
on a fixed smooth
 manifold $M$. The reason we have chosen this name in favor of the alternative, stems from the
  fact that our ``stacks" are in fact generalizations of the  above familiar 1-categorical 
concept to $(\infty,1)$-categories. Given a 1-truncated stack (in our sense), if one restricts to 
covers by convex subsets of a smooth manifold $M\in \Sh(\cartsp)$, one recovers the
 usual notion of a smooth stack on $M$. We have also chosen stack, instead of sheaf, since our 
 objects do not satisfy the \emph{strict} gluing condition, but only satisfy descent in general.

\medskip
In Section \ref{Sec vb} we will first review some of the theory and techniques from derived geometry and then 
use these to study bundles and their classifications in stacks. 
To help dissolve some of these conceptual difficulties, we provide some intuition as to how to bridge the gap between 
familiar smooth objects and smooth stacks. Let us consider the familiar example of the category of smooth manifolds 
$\mathscr{M}{\rm an}$. It is well known that this category is very poorly behaved under category theoretic constructions
(see \cite{Stacey} for a detailed discussion). For example, pullbacks and pushouts of smooth maps between manifolds 
are very seldom smooth manifolds themselves. 
However, there is a fully faithful embedding $\mathscr{M}{\rm an}\into \Sh(\cartsp)$ into the category of smooth sheaves. 
In contrast, the latter category is very well-behaved under these constructions; in fact, it is 
complete and cocomplete (i.e. has all limits
and colimits). Thus we see 
that, by identifying the category of manifolds with a subcategory of sheaves, we can perform categorical constructions 
and identify the results with a \emph{smooth sheaf}. It is worth emphasizing that the fact that the category of smooth 
manifolds embeds fully faithfully into sheaves means that any construction which uses only smooth manifolds and maps 
between these holds equally well in the category of smooth sheaves.

\medskip
The category of smooth stacks $\Sh(\cartsp)$ is a higher categorical analogue of the category of smooth sheaves. Thus, we can view smooth sheaves as being a subcategory of smooth stacks in the same way that sets, viewed as discrete topological spaces, are a subcategory of topological spaces. So we see that smooth manifolds can really be thought of as certain ``discrete" smooth stacks. At first glance this might be confusing since, as a space, a smooth manifold is far from 
being discrete and it may have highly nontrivial higher homotopy groups. 
\footnote{ \label{foo 1} Note that this is analogous to the classical situation of a free action of 
a finite group $\Gamma$ on a space $M$, where $\pi_i(M)=\pi_i(M/ \Gamma)$ for $i \geq 2$
while $\pi_1(M/\Gamma) \hookrightarrow \pi_1(M)$ is an injection. Indeed, 
with $\Gamma$ viewed as a discrete space, the fibration 
${\cal F} \to {\bf M}_{\cal F} \to {\bf M}$ is similar to  
$\Gamma \to M/\Gamma \to M$ in that the only effect of the fiber is on $\pi_0$. 
Therefore, a sheaf in the world of stacks is, in the above sense, the analog of 
a discrete space.} Here we come to an interesting phenomenon 
that occurs in smooth stacks. The  higher homotopy data is absorbed by the geometry of the manifold. More precisely, 
smooth stacks have a way of separating out the homotopy theoretic data arising via the smooth structure from that which 
came from more combinatorial constructions.  The smooth data is absorbed by the simplicially discrete
\footnote{Note again that being discrete in spaces and in stacks have totally different meanings. 
See the discussion preceding footnote \ref{foo 1}.}
sheaves, while the combinatorial data is encoded in the simplicial homotopy groups of the stack. So, in a sense that can be 
 made precise, a smooth stack is a level-wise  smooth object whose levels are glued together according to simplicial rules.
  As an example of a smooth stack with nontrivial combinatorial data, the smooth stack   $\BB^{n-1} {\rm U}(1)$, which is a smooth 
  version of the Eilenberg-MacLane space $K(\ZZ,n)$, has nontrivial simplicial homotopy group in degree $n$
  (see \cite{FSSt} \cite{SSS3}).

\medskip
The Madsen-Tillman construction \cite{MT} involves an interplay between the tangent bundle 
$TM$ and the normal bundle $T^{\perp}M$ of an embedding $j: M \hookrightarrow \R^{d+N}$ of a manifold $M$ (of dimension $d$) into 
a Euclidean spaces of large dimension ($N \to \infty$). The two bundles combine together to provide 
a decomposition of vector bundles $TM \oplus T^{\perp}M\cong j^* T\R^{d+N}$, exhibiting the two original bundles 
as bundles of vector spaces which are orthogonal complements in a large vector space. 
For an abstract manifold, these relationships implicitly involve a number of choices. In particular, we have a cover of $M$ by coordinate patches and a trivialization of both bundles over the local patches. Smooth stacks provide a way for making all this information manifest. At first glance, keeping track of such data might seem to be irrelevant, if not pedantic. Nevertheless, these choices turn out to be crucial in the Pontrjagin-Thom construction in smooth stacks. The reason essentially boils down to the following observation. If one is given a \emph{submanifold} $M\subset \RR^{d+N}$ of Euclidean space, there is a canonical map to the Grassmannian classifying the tangent bundle. This map has nothing to do with a choice of cover: it simply maps a point in the manifold to the tangent space (viewed as a subspace of $\RR^{d+N}$). If we view $M$ as an abstract manifold, there is no canonical choice of map and the data of local trivializing patches needs to be taken into account. The key observation in the proof of our main theorem is that if we unravel the local trivializing information, then the inverse to the Pontryagin-Thom map can be calculated locally and glued together to give a submanifold. Moreover, the local problem turns out to be nothing more than of finding a regular value of a smooth map $f:D^{N+d}\to D^N$ and witnessing the local inverse as the corresponding Pontryagin submanifold of $D^{d+N}$.

\medskip 
In Section \ref{Sec gras} we develop the classification of vectors bundles and their complements
in stacks. Hence we consider stacky generalizations of the classifying space as $\mathbf{B}{\rm O}(n)$,
of the universal bundle as $\mathbfcal{U}(d,N)$, of its orthogonal complement as $\mathbfcal{U}^\perp(d,N)$, of Grassmannians as $\mathbf{Gr}(d, N)$, and of Stiefel manifolds as $\mathbf{V}(d, N)$. 
The main point will be that, in order to compute the inverse map, we need various models for the corresponding smooth manifolds in stacks. With the exception of $\BB{\rm O}(d)$, these stacks will all be equivalent to their classical counterparts. \footnote{Here we mean, equivalent in the category of smooth stacks.} We hope that these constructions will be of use in their own right. 
These are summarized in Table 1, which we hope would also be useful for notation. 


\medskip
\medskip
{
\begin{table} 
\hspace{5mm}
\small
\begin{tabular}{ |>{\centering\arraybackslash}m{7cm}||>{\centering\arraybackslash}m{7.5cm}
|}
\hline
\begin{center}
{\bf Stack} 
\end{center} &  \begin{center} {\bf Description} \end{center}
\\
\hline
\hline
Classifying stack of orthogonal bundles   
${\bf B} {\rm O}(n)$ & Homotopy orbit stack of  the smooth \newline 
sheaf ${\rm O}(n)$ acting on the point 
\\
\hline
Universal bundle \newline
$\mathbfcal{U}(d)\to \BB{\rm O}(d)$& $\mathbfcal{U}(d)\simeq \RR^d/\!/{\rm O}(d)$  
\\
\hline
Real Grassmannian stack
$\mathbf{Gr}(d,N)$  & 
$
\begin{aligned} 
\mathbf{Gr}(d,N)
&:=
{\rm O}(d+N)/\!/{\rm O}(d)\times {\rm O}(N)
\\
&\simeq 
{\rm O}(d+N)/{\rm O}(d)\times {\rm O}(N)
=:{\rm Gr}(d, N)
\end{aligned} 
$
\\
\hline
Real Stiefel stack 
$\mathbf{V}(d,N)$  
& 
$
\begin{aligned}
\mathbf{V}(d,N) &:= 
{\rm O}(N)/\!/{\rm O}(N-d)
\\
& \simeq
{\rm O}(N)/{\rm O}(N-d)=:
{\rm V}(d,N)
\end{aligned}
$
\\
\hline
Universal vector bundle \newline $\mathbfcal{U}(d,N)\to {\bf Gr}(d,N)$ 
&
$
\big(\RR^d\times {\rm V}(d,N)\big)/\!/{\rm O}(d)
$
\\
\hline
Universal complement bundle 
$\mathbfcal{U}^{\perp}(d,N)\to {\rm Gr}(d,N)$ 
&
$
\big((\RR^d)^{\perp}\times {\rm V}(N,d)\big)/\!/{\rm O}(N)
$
\\
\hline
{\rm Thom stack} $\mathbf {Th}(\eta)$ of a vector bundle $\eta\to X$ 
&
$
\begin{aligned}
 \mathbf {Th}(\eta) &:=D(\eta)/\!/S(\eta)
 \\
&\simeq D(\eta)/S(\eta) =:{\rm Th}(\eta)
 \end{aligned}
 $
\\
\hline
Stacky sphere \newline $D^n/\partial D^n$
 & 
 $D(V)/S(V)$, ~~${\rm dim}(V)=n$ 
 \\
\hline
Thom stack of the universal bundle 
 $\mathbfcal{U}(d)\to \BB{\rm O}(d)$ 
 &
 ${\rm Th}(\mathbfcal{U}(d)) \simeq  (D^d/\partial D^d)/\!/{\rm O}(d)$
 \\
 \hline
 Thom stack of the universal bundle 
$\mathbfcal{U}(d,N)\to {\rm Gr}(d,N)$ 
&
${\rm Th}(\mathbfcal{U}(d,N))\simeq 
  \big(D^d/\partial D^d\wedge {\rm V}(d,N)\big)/\!/{\rm O}(d)$
\\
\hline
Thom stack of the universal bundle 
$\ \mathbfcal{U}^{\perp}(d,N)\to {\rm Gr}(d,N)$ 
& ${\rm Th}(\mathbfcal{U}^{\perp}(d,N))\simeq  
\big(D^N/\partial D^N\wedge {\rm V}(N,d)\big)/\!/{\rm O}(N)$
\\
\hline
Stacky Madsen-Tillman spectrum $ \MM \TT(d) $
&
$
\begin{aligned}
 \MM \TT(d) & :={\rm Th}(\mathbfcal{U}^{\perp}(d,N))
\\ 
 & \simeq \big(D^N/\partial D^N\wedge {\rm V}(N,d)\big)/\!/{\rm O}(N)
 \end{aligned}
 $
  \newline
\\
\hline
Smooth motivic model for the Thom spectrum 
$\mathbf{M}{\rm O}$
&
$\mathbf{M}{\rm O}:={\rm Th}(\mathbfcal{U}(d))\simeq
 (D^d/\partial D^d)/\!/{\rm O}(d)$ 
\\
\hline
$D^1/\partial D^1$-suspension spectrum  
$\Sigma^{\infty}_{{}_{D^1/\partial D_1}}X$ 
& 
Successively smashing with the stacky circle 
$D^1/\partial D^1$
\\
\hline
 $D^1/\partial D^1$-infinite loop stack 
 $
\Omega^{\infty}_{{}_{D^1/\partial D_1}}X(n)
$
 &
 \begin{center}
$
\Omega^{\infty}_{{}_{D^1/\partial D_1}}X(n)
:=
\colim_{n}\mathbf{Map}_{+}\big(D^n/\partial D^n,X(n)\big)
$ 
\end{center}
\\
\hline
Smooth concordance
 category  
 with \newline
 ordering on $[0,1]$: ${\rm Conc}^{>}(X)$ 
 &
 Diagram on the full substack 
 of collared \newline
 and ordered maps
\\
\hline
 $\BB{\rm Conc}^{>}(X)$
 &
 Realization of the nerve of groupoid \newline
${\rm Conc}^{>}(X)$
\\
\hline
Smooth cobordism category \newline
$\mathscr{C}{\rm ob}_d$ 
&
Quotient of of embeddings and diffeomorphisms
 \\
\hline
Smooth cobordism stack \newline 
$\BB\mathscr{C}{\rm ob}_d$ 
&
Stacky realization of  $\mathscr{C}{\rm ob}_d$ 
\\
\hline
Smooth stack of bundle maps 
${\bf Bun}(T^{\perp}M,\mathbfcal{U}(d,N))$
&
Stackification of $\map(T^{\perp}M, \mathbfcal{U}(d,N))$ 
on those maps which are bundle maps
\\
\hline
Smooth stack of embeddings with local trivializations 
$\mathbf{Emb}_{\rm loc}(M,\RR^{d+N-1})$ &  
Homotopy pullback
of \newline
 $\mathbf{Bun}\big(T^{\perp}M,\mathbfcal{U}^{\perp}(d,N)\big)$
and 
${\rm Emb}(M,\RR^{d+N-1})$
\\
\hline
\end{tabular}
\caption{Various stacks and categories used in the paper. As we show, 
most of these stacks end up being zero-truncated.}
\end{table}
}

\medskip
The techniques developed by Galatius, Madsen, Tillman and Weiss 
\cite{MT}\cite{MW}\cite{GMTW}, 
which led to a proof \cite{MW} of the Mumford conjecture \cite{Mum}, have seen a number of fascinating applications in recent years (see \cite{Ralph} for an exposition). Notably, their work has led to a deeper understanding of topological quantum field theories (TQFT's), which has since been far 
generalized by the work of Lurie on the cobordism hypothesis \cite{Lur2} (see \cite{Fr1} for an exposition). 
One facet, which has not been so extensively studied in the literature, is the application of their techniques to 
quantum field theories (QFT's) which have additional geometric data (such as connections, Riemannian metrics, ect.). These types of QFT's are inherently not topological, but at the same time are really the field theories which one needs to understand in 
physical applications. Nevertheless, there has been some work in this direction, for example 
Cohen, Galatius and Kitchloo \cite{CGK} have applied the these techniques to moduli spaces of flat 
connections and Ayala \cite{A} has used sheaf theoretic techniques to study some of the effects of 
adding other geometric data to the mix. However, the end result in both cases is an identification of the \emph{geometric realization} of certain cobordism categories with a \emph{space}. 
As such, the homotopy type of these realizations often lose a large amount of the geometric data.
We plan for a future study of field theories in our current context. 
\footnote{We thank Ralph Cohen for very useful discussions in that direction.}

 \medskip
Adding connections to the moduli spaces, when viewed from a topological point of view
might seem like nothing has happened, as the space of connections is affine, i.e. 
contractible. However, we certainly would like for the connections to have a nontrivial 
effect. The way to improve upon this situation is to describe these categories not just as topological 
categories, but a smooth sheaves of topological categories, or smooth stacks. For our purposes, it 
will be better to work with \emph{simplicial sheaves} as a model for smooth stacks, rather than 
topological valued sheaves and we will make the transition back to topological spaces by 
geometrically realizing.

\medskip
As an example of how this occurs, consider the sheaf of Lie-algebra valued 1-forms 
$\Omega^1(-;\mathfrak{g})$. This sheaf encodes globally defined connections on 
$G$-principal bundles: Given a smooth manifold $M$, the sheaf condition provides 
us with an identification
$$
\pi_0\map(M,\Omega^1(-;\mathfrak{g}))\simeq \Omega^1(M;\mathfrak{g})\;.
$$
This example highlights that homotopy theory of smooth stacks captures the geometric 
data encoded by global connections. However, under the geometric realization functor
$$
\vert{\Pi}\vert: {\Sh}(\cartsp)\longrightarrow 
\sset
\longrightarrow
 \mathscr{T}{\rm op}\;,
$$
the space $\vert \Pi \Omega^1(-;\mathfrak{g})\vert$ can be identified with the affine space of $\mathfrak{g}$-valued 1-forms on  $\RR^{\infty}$ (understood as the colimit of $1$-forms 
on $\RR^N$ as $N\to \infty$). Being a vector space, this space is contractible and therefore the space of maps $\map(M,\vert \Pi \Omega^1(-;\mathfrak{g})\vert)$ is contractible. Thus, we see that the homotopy type of the realization loses the geometric information encoded by the differential forms.

\medskip

Starting in Section \ref{Sec GMTW}, we refine the classifying space of the cobordism category to smooth stacks. It turns out that in the category of smooth stacks, the corresponding classifying stack is zero-truncated and is therefore equivalent to its sheaf of connected components. This will allow us to only define the abstract Pontryagin-Thom collapse map at only at the level of connected components and we can mimick the construction in spaces to produce the map.

\medskip
In Section \ref{Sec PT} we provide a smooth refinement of the 
Pontrjagin-Thom construction to the category of smooth stacks. Upon geometrically 
realizing our  stacks, we recover the usual collapse map. 
This map will be a morphism of sheaves \eqref{abstract collapse map}
$$
\xymatrix{
{\bf PT}:\widetilde{\pi}_0\big(\BB\mathscr{C}{\rm ob}_{d}\big)
\; \ar[r] & \;
\widetilde{\pi}_0\big(\BB{\rm Conc}^{<}\big(\Omega^{\infty-1}_{{}_{D^1/\partial D^1}} {\MM\TT}(d)\big)\big)
}\;,
$$
where the source stack is the stacky analogue of the classifying space of the cobordism category defined in Section \ref{Sec GMTW} 
(see Definition \ref{scc}) and the target category is defined in 
Section \ref{Sec mot} (see Corollary \ref{realization of concordance}).

\medskip

\begin{introthm}{Theorem}[\bf{Stacky Pontrjagin-Thom equivalence}]
The map ${\bf PT}$ induces a weak equivalence of smooth stacks
$$
\xymatrix{
{\bf PT}:\BB\mathscr{C}{\rm ob}_{d}
\; \ar[r]^-{\simeq} & \;
\BB{\rm Conc}^{>}\big(\Omega^{\infty-1}_{{}_{D^1/\partial D^1}}{\MM\TT}(d)\big)
}\;.
$$

\end{introthm}

This is Theorem \ref{theorem 1} in Section \ref{Sec PT}. 
Consequently, the result we will prove in Proposition \ref{Prop real} leads to a weak equivalence 
$\vert \BB \mathscr{C}{\rm ob}_{\theta}\vert \to  {\rm B} \vert \mathscr{C}{\rm ob}_{\theta}\vert 
\to {\rm B}\mathscr{C}{\rm ob}_{\theta}$. Similarly, Proposition \ref{geometric realization of spectrum} 
gives a weak equivalence $\vert \BB{\rm Conc}^{>}(\Omega^{\infty-1}\MM \TT(\theta)) \vert\to \Omega^{\infty-1}{\rm MT}(\theta)$. By the theorem, ${\bf PT}$ 
defines a weak equivalence of smooth stacks and since the geometric realization functor sends weak 
equivalences 
to weak equivalences we have the following.

\begin{introthm}{Corollary}[\bf{Galatius-Madsen-Tillman-Weiss}]
The map ${\bf PT}$ induces a weak equivalence
$$
\xymatrix{
{\rm PT} : {\rm B}\mathscr{C}{\rm ob}_{\theta}
\; \ar[r]^-{\simeq}& \;
\Omega^{\infty-1} {\rm MT}(\theta)
}\;.
$$
\end{introthm}

For a tangential structure $\theta$, which is induced by a faithful representation $\theta:G\into {\rm O}(d)$, we have a generalization of the first theorem.
The following is Theorem \ref{theorem 2} in Section \ref{Sec G}.

\begin{introthm}{Theorem}[{\bf Stacky Pontrjagin-Thom equivalence with $\theta$-structure}]
The map  
$$
\xymatrix{
{\bf PT}_{\theta}:\BB\mathscr{C}{\rm ob}_{\theta}
\; \ar[r] & \;
\BB{\rm Conc}^{>}\big(\Omega^{\infty-1}_{{}_{D^1/\partial D^1}}{\MM\TT}(\theta)\big)
}
$$
is a weak equivalence of smooth stacks. 
\end{introthm}

Going beyond $G$-structure, we can also add more refined geometric structure in our setting 
(Section \ref{Sec geom}). We can consider, for example, the cobordism category with objects 
smooth manifolds whose tangent bundles are equipped with a connection and whose morphisms
 are bordisms with tangent bundles equipped with connections extending those of the bounding 
 manifolds. We can also consider manifolds equipped with 
Riemannian structure, symplectic structure, complex structure and so on.

\begin{introthm}{Theorem}[{\bf Stacky Pontrjagin-Thom equivalence with extra geometric structure}]
The map
$$
\xymatrix{
{\bf PT}^{\cal F}:\BB\mathscr{C}{\rm ob}_{d}^{{\cal F}} 
\; \ar[r] & \;
 \BB{\rm Conc}^{>}\big(\Omega_{{}_{D^1/\partial D^1}}^{\infty-1}{\MM\TT}(d)_{\cal F}\big)
 }
$$
is a weak equivalence of smooth stacks.
\end{introthm}

This is Theorem \ref{final thm} in Section \ref{Sec geom}. 
In \cite{A}, the topological category of cobordisms with ${\cal F}$-structure (where ${\cal F}$ 
is a sheaf on the site of smooth manifolds with values in $\mathscr{T}{\rm op}$) is defined. 
This category has space of objects equivalent to the coproduct of homotopy orbit spaces 
${\cal F}(M)/\!/{\rm Diff}(M)$ where $M$ ranges through diffeomorphism classes of 
manifolds of fixed dimension. Similarly, the morphisms are identified with coproducts of 
${\cal F}(W)/\!/{\rm Diff}(W)$, where $W$ is a bordism between manifolds. 
Ayala proves that there is a weak homotopy equivalence
$$
{\rm B}\mathscr{C}{\rm ob}_{d}^{{\cal F}}\simeq \Omega^{\infty-1}{\rm MT}(d)_{\cal F}\;,
$$
where ${\rm MT}(d)_{\cal F}$ is the spectrum which at level $n$ is given by 
${\rm Th}(p^*_{\cal F}\;{\cal U}^{\perp}(d,N))$ with the map
$p_{\cal F}: {\rm Gr}(d,N)_{{\cal F}(\RR^d)}\to {\rm Gr}(d,N)$, where
${\rm Gr}(d,N)_{{\cal F}(\RR^d)}$ is defined analogous to our definition in smooth 
stacks, but in the category of spaces. 
Theorem \ref{final thm} then achieves a similar goal as Ayala upon geometric realization.

%

\medskip
We can also consider mixed situations, i.e. when we have both a tangential structure, as a 
$\theta$-structure, and geometric data, as an $\cal{F}$-structure. 
\begin{itemize}

\vspace{-2mm}
\item These structures might
a priori not be directly related, such as a Riemannian metric and a Spin structure. In this 
case we would have an action of the Spin group on the space 
${\rm Riem}(n) \times \Omega^1(-; \mathfrak{so}(n))$ via the surjection 
${\rm Spin}(n) \to {\rm O}(n)$. 

\vspace{-2mm}
\item We could also consider
situations when the $\theta$-structure and the $\cal{F}$-structure are required to be 
compatible. For instance, we could consider a Spin structure together with a 
compatible connection. 
\end{itemize}

In general, we describe such combinations via $(\theta, {\cal F})$-structures, in which 
case we have yet another generalization.

\begin{introthm}{Theorem}[\bf{Stacky Pontrjagin-Thom equivalence with geometric and $\theta$-structures}]
\label{Thm most}
The map
$$
\xymatrix{
{\bf PT}_{\theta}^{{\cal F}}:\BB\mathscr{C}{\rm ob}_{\theta}^{{\cal F}} 
\; \ar[r] &  \; 
\BB{\rm Conc}^{>}\big(\Omega_{{}_{D^1/\partial D^1}}^{\infty-1}{\MM\TT}(\theta)_{ {\cal F}}\big)
}
$$
is a weak equivalence of smooth stacks.
\end{introthm}

%


In \cite{RS}, a different notion of parametrized cobordism is introduced,  
where the parametrization is with respect to the $\theta$-structure, while 
ours is more about geometrically refining  the cobordism category 
in \cite{GMTW}.

\vspace{.5cm}
\par{\large \bf Acknowledgement.} 
The authors would like to thank Ralph Cohen for very useful discussions and encouragement. The authors had been working on differential refinement of cobordism, as a generalized cohomology theory, when they attended the 
Weiss 60 meeting in Lisbon (H.S.) and the Young Topologists meeting in 
Copenhagen (D.G.) in summer 2016, which inspired them to refine the 
Galatius-Madsen-Tillman-Weiss setting. The authors thank the organizers, 
speakers, and participants of the two meetings. 
The authors are grateful to Oscar Randal-Williams and Urs Schreiber for 
a careful reading of the first draft and for useful suggestions and also to 
David Ayala for useful discussions about his work.

\section{Vector bundles and their classification in derived geometry} 
\label{Sec vb}

\subsection{Recollection and techniques from smooth stacks}
\label{Sec recol}

In this section, we provide a quick review of smooth stacks and the basic properties used frequently throughout 
this paper. The idea of the proofs of the main theorems are more or less easy to understand, at least conceptually, 
without having to fully grasp all the technical details involved. Thus, the reader who is not interested in the technical 
details of the proofs might wish to only have a quick glance at this section.

\medskip
In essence, smooth stacks provide us with a convenient place to do parametrized homotopy theory on 
objects which glue with respect to some local data (e.g. manifolds, algebraic varieties, etc.). To begin, 
we start with a category equipped with the notion of covering families (i.e. a \emph{site}, 
see for example \cite{MM} 
for details). This category provides us with our parametrizing spaces and the coverings tell us how 
to glue along local data. Since we are concerned with \emph{smooth} stacks, we will restrict our 
attention to prestacks on the following site.
 \begin{definition}
 [The category of Cartesian spaces and the category of smooth prestacks]
{\bf (i)} Let $\cartsp$ be the site with objects the convex subspaces of Euclidean space and morphisms 
the smooth maps between them. \footnote{By smooth map, we mean in the sense of manifolds with 
corners. Since all our objects are subsets of Euclidean spaces, this is equivalent to requiring the 
existence of a smooth extension  on an open neighborhood.} We topologize $\cartsp$ by taking good 
open covers  (i.e. covers with contractible finite intersections)  of convex subspaces. 

\noindent{\bf (ii)} Let $\mathscr{P}{\rm sh}_{\infty}(\cartsp):=[\cartsp^{\rm op},\sset]$
 denote the category of contravariant functors
$$
X:\cartsp^{\rm op}\longrightarrow \sset\;,
$$
with natural transformations between them. We call this category, the {\rm category of smooth prestacks}.
\end{definition}

\medskip
Essentially all the tools that one has in the category of simplicial sets can be applied to prestacks. For instance, 
we have the following (see \cite{Jar}). 
\begin{definition} 
[Presheaf of connected components]
The \emph{presheaf of connected components of a prestack} $X$ is defined 
as the presheaf \footnote{As in the case of topological spaces, $\widetilde{\pi}_0$ 
does not admit a group structure in general.}
$$
\pi_0(X):\cartsp^{\rm op}\longrightarrow \set\;,
$$
which assigns each convex subspace $U\in \cartsp$ to the set of simplicial connected components 
$\pi_0(X(U))$. 
\end{definition} 

Note that $\pi_0$ is functorial and $X$ is functorial in $U$, which indeed defines 
a presheaf. In a similar way, we have the following (see \cite{Jar}). 

\begin{definition}
[Presheaf of homotopy groups]
Given a basepoint of a prestack $x:\ast\to X$,  
the \emph{presheaf of homotopy groups} can be defined, for $m>0$, 
$$
\pi_m(X,x):\cartsp^{\rm op}\longrightarrow \set\;,
$$
as the functor which assigns to each $U\in \cartsp$ the simplicial 
homotopy group $\pi_m(X(U),x_{U})$. 
\end{definition} 

Here the basepoint $x_{U}$ is prescribed by the natural transformation $x:\ast\to X$. Given that these are presheaves on the site $\cartsp$, equipped with the notion of covering families (i.e. the good open covers 
of convex subspaces), we can \emph{sheafify} the resulting presheaves with respect to the coverage 
(see \cite{Jar} for details).

\medskip
These observations lead us to a natural definition for a model structure on prestacks. With the injective model structure, this was first defined by Jardine in \cite{Jar}. The projective model structure is considered in 
other places (for instance, in \cite{DI}\cite{FSSt}\cite{Urs}). For us, it will turn out that the projective model structure is slightly more convenient for 
calculations and hence we will adopt it as our model structure.

\begin{definition}
[The category of smooth stacks]
\label{smooth stacks}
We define the $\infty$-{\rm category of smooth stacks} 
$$
{\Sh}(\cartsp):={\rm fib}([\cartsp,\sset]^{\rm proj, loc})
$$ 
as the full subcategory on fibrant 
objects of the smooth prestack category $[\cartsp,\sset]^{\rm proj, loc}$, equipped with the 
projective model structure on functors (i.e. the fibrations are objectwise Kan fibrations) 
and localized at the morphisms of prestacks $X\to Y$ which induce isomorphisms on all 
sheaves of homotopy groups.
\footnote{If the stack has more than one connected component, this means that the map must induce an 
isomorphism on sheaves of homotopy groups, with basepoint taken in each connected component.}
\end{definition} 

We will often need mappings between stacks. 

\begin{definition} 
[Mapping spaces of stacks] The mapping spaces between two smooth stacks is given via the simplicial model 
structure on the localization. More precisely, given two smooth stacks $X$ and $Y$, the mapping 
space is defined as $\map(X,Y):=\map(Q(X),Y)$, where $Q(X)$ denotes some cofibrant 
replacement of $X$.
\end{definition}

Later we will also encounter mapping stacks which are not spaces (see Proposition 
\ref{eor pro} and Definition \ref{susloop}).

\medskip
Note that if a map $X\to Y$ induces an isomorphism at the level of \emph{presheaves} of homotopy groups, this immediately 
implies that the corresponding sheafifications are isomorphic. The converse is not necessarily true, however. Thus, the class 
of weak equivalence in ${\Sh}(\cartsp)$ is strictly larger than the objectwise weak equivalences, i.e. an equivalence 
of simplicial sets for each $U$.

\begin{remark}
[Other models for smooth stacks]
There are many possible presentations for the $\infty$-category defined above. For example, one could take localization hammocks on the presheaf category with local weak equivalences (see \cite{DHKS}). One could also consider smooth stacks on the site of smooth manifolds, which would lead to an equivalent $\infty$-category (see \cite{Urs}). 
We have chosen the above model since it seems to simplify some calculations. 
\end{remark}

At first glance, it might seem awkward that we have defined the category of smooth stacks as the fibrant objects with respect to some model structure. Indeed, if we are to follow traditional literature, we would like for a stack to be characterized by some sort of descent condition. 
Thus it would seem  that we have departed from what has been traditionally done in the literature, namely characterize a stack via descent condition. However, happily the two notions agree. Indeed, Dugger, Hollander and Isakson show in \cite{DI}, that
the fibrant objects are precisely the objectwise Kan complexes which satisfy descent with respect to 
all hypercovers. Put another way this means that, for a functor $F:\cartsp^{\rm op}\to \sset$
which takes values in Kan complexes, we have an equivalence
$$
F(X)\simeq {\rm holim}\left\{
\xymatrix{
\hdots & \ar@<.15cm>[l] \ar@<-.15cm>[l] \ar@<.05cm>[l]\ar@<-.05cm>[l]   F(Y_2) & \ar@<.1cm>[l] \ar[l] \ar@<-.1cm>[l]   F(Y_1) &  \ar@<-.05cm>[l] \ar@<.05cm>[l]   F(Y_0)
}\right\}\;,
$$
for all hypercovers $Y_{\bullet}\to X$, if and only if $F$ is an object in ${\Sh}(\cartsp)$. 
Using the basic properties of the mapping space in a simplicial model category, it is easy to show 
(see again \cite{DI}) that this is equivalent to saying that $F$ is fibrant in Bousfield localization 
taken at the hypercovers $Y_{\bullet}\to U$. 

\begin{remark}
[Hypercover vs. {\v C}ech nerve]
In general, hypercovers form a strictly larger class than the ordinary {\v C}ech nerves of covers. 
However, over the site of Cartesian spaces, the fiber products of cover inclusions are again 
Cartesian spaces. Thus, every hypercover is actually a {\v C}ech nerve and we can safely 
restrict to this case.
\end{remark}

Notice that if we fix a functorial fibrant replacement functor $L$ on the 
smooth prestack model category
 $[\cartsp^{\rm op},\sset]^{{\rm proj},{\rm loc}}$, such a functor gives us a way to turn 
 a prestack into a smooth stack. Clearly $L^2\simeq L$ and $L$ defines a left $\infty$-adjoint 
 to the inclusion. Moreover, it can be shown that $L$ is also left exact (i.e. preserves finite limits) 
 \cite{Lur}. Summarizing, there is an $\infty$-adjunction between the category of stacks and 
 that of prestacks
$$
\xymatrix{
\Sh(\cartsp)\  \; \ar@<-.15cm>@{^{(}->}[r]_-{i} & \; \ar@<-.15cm>[l]_-{L}\PSh(\cartsp)\;,
}
$$
with $L$ preserving finite  $\infty$-limits. Any such functor $L$ is called 
a \emph{stackification functor}. This functor is not unique, but is unique up to canonical equivalence.

\begin{example}
[Smooth manifold as a smooth stack]
For a smooth manifold $M$ (possibly with boundary, or even corners), let $C^{\infty}(-,M)$ be the 
smooth sheaf given by sending a convex space $U$ to $C^{\infty}(U,M)$, where smoothness is 
defined as the usual notion of a smooth map between manifolds with corners. By extending this sheaf 
to be constant in the simplicial direction, we get an object of ${\Sh}(\cartsp)$. This assignment 
gives rise to a fully faithful embedding
$$
\xymatrix{
\mathscr{M}{\rm an} \; \ar@{^{(}->}[r] & {\Sh}(\cartsp)
}.
$$
Henceforth, as we mentioned earlier, whenever we speak of a smooth manifold
 we will mean the corresponding object in ${\Sh}(\cartsp)$. 
\end{example}

For a general smooth prestack, and not just an object in the site $\cartsp$, there is a more general 
notion of covering which reduces to the usual notion of {\v C}ech nerve when restricted to convex open subsets of $\R^n$. 
A \emph{local epimorphism} is a morphism of smooth prestacks $f:X\to Y$, such that for each $U\in \cartsp$ and each 
map $i:U\to X$, the iterated fiber product 
$$
\xymatrix{
\hdots  \ar@<.15cm>[r] \ar@<.05cm>[r] \ar@<-.05cm>[r] \ar@<-.15cm>[r]  & X\times_{Y} X\times_{Y} X\ar@<.1cm>[r]\ar[r]\ar@<-.1cm>[r] & X\times_Y\times X \ar@<-.05cm>[r] \ar@<.05cm>[r] & X \ar[r]^-{f} &Y
\\
\hdots  \ar@<.15cm>[r] \ar@<.05cm>[r] \ar@<-.05cm>[r] \ar@<-.15cm>[r]  & \coprod_{\alpha\beta\gamma}U_{\alpha\beta\gamma}  \ar[u] \ar@<.1cm>[r]\ar[r]\ar@<-.1cm>[r] & \coprod_{\alpha}U_{\alpha\beta} \ar@<-.05cm>[r] \ar@<.05cm>[r] \ar[u]& \coprod_{\alpha}U_{\alpha} \ar[u]\ar[r] & U\ar[u]_-{i}
\;,
}
$$
is a {\v C}ech nerve of a cover $\coprod_{\alpha}U_{\alpha}\to U$. The following example illustrates how useful this more general notion can be when working over the small site of Cartesian spaces.

\begin{example}
[Cover for smooth manifolds]
Let $M$ be a smooth manifold and let $\{U_{\alpha}\}$ be an open cover of $M$ with contractible finite intersections (i.e. a good open cover). Then we can form the {\v C}ech nerve of the cover and regard it as a simplicial diagram in prestacks,
$$
\xymatrix{
\hdots \ar@<.1cm>[r] \ar[r] \ar@<-.1cm>[r] &
\coprod_{\alpha\beta}U_{\alpha\beta}  \ar@<-.05cm>[r] \ar@<.05cm>[r] & 
\coprod_{\alpha}U_{\alpha} \ar[r]& M
}\;.
$$
Since $M$ itself is not representable, this {\v C}ech nerve cannot define a covering. However, it is 
clear that this object defines a local epimorphism of prestacks, since given a smooth map $f:U\to M$, 
the iterated fiber products is just the nerve of the cover $\coprod_{\alpha}f^{-1}(U_{\alpha})$.
\end{example}

Every local epimorphism $f:X\to Y$ comes equipped with a map
$$
{\rm hocolim}_{\Delta^{\rm op}}\check{C}(f)\longrightarrow Y\;,
$$
and in the category of prestacks there is not too much one can say about this map. However, as a consequence 
of the charaterization of fibrant objects via descent, one sees that fibrantly replacing (or stackifying) this map
gives an equivalence of smooth stacks. On the other hand, the stackification functor $L$ is both left exact and a left adjoint (in the $(\infty,1)$-sense). As a consequence, we have an equivalence
$$
{\rm hocolim}_{\Delta^{\rm op}}\check{C}(L(f))\overset{\simeq}{\longrightarrow} L(Y)\;.
$$
The morphisms of \emph{stacks} $f:X\to Y$, which induces an equivalences out of their {\v C}ech nerves 
are called \emph{effective epimorphisms}. The above argument shows that effective epimorphisms are 
precisely the stackification of local epimorphisms.

\begin{example}
[Cover for homotopy quotients]
Let $G$ be a Lie group acting on a smooth stack $X$. Let $q:X\to X/\!/G$ be the canonical map to the homotopy quotient by the $G$-action. The {\v C}ech nerve of the quotient takes the form
$$\xymatrix{
\hdots  \ar@<.15cm>[r] \ar@<.05cm>[r] \ar@<-.05cm>[r] \ar@<-.15cm>[r] & X\times G\times G\ar@<.1cm>[r]\ar[r]\ar@<-.1cm>[r] & X\times G \ar@<-.05cm>[r] \ar@<.05cm>[r] & X\ar[r] & X/\!/G
}\;.
$$
There is an internal characterization of $(\infty,1)$-topos due to 
Rezk \cite{Re} and Lurie \cite{Lur} and one of the characterizing axioms asserts that quotients by $G$-actions are always effective. In other words, the above diagram is a homotopy colimiting diagram. Thus, the morphism $q$ in prestacks must be a local epimorphism and can be thought of as a cover for the homotopy quotient.
\end{example}

In the above example, if $X$ is a smooth \emph{sheaf} and $G$ acts freely on $X$, then the homotopy quotient $X/\!/G$ is modeled by the strict quotient sheaf $X/G$. So we arrive at the identification 
$$
{\rm hocolim}_{\Delta{\rm op}}\left\{\xymatrix{
\hdots  \ar@<.15cm>[r] \ar@<.05cm>[r] \ar@<-.05cm>[r] \ar@<-.15cm>[r] & X\times G\times G\ar@<.1cm>[r]\ar[r]\ar@<-.1cm>[r] & X\times G \ar@<-.05cm>[r] \ar@<.05cm>[r] & X
}
\right\}
\simeq X/G
$$
On the other hand the homotopy colimit on the right can be calculated by stackifying the objectwise homotopy colimit. This, in turn, can be computed via the diagonal of the relevant bisimplicial complex. In this case, both $G$ and $X$ are zero-truncated and the bisimplicial complex is completely degenerate in one direction. Applying the 
stackification functor $L$ gives the identification 
\(\label{shdc}
L\left\{\xymatrix{
\hdots  \ar@<.15cm>[r] \ar@<.05cm>[r] \ar@<-.05cm>[r] \ar@<-.15cm>[r] & X\times G\times G\ar@<.1cm>[r]\ar[r]\ar@<-.1cm>[r] & X\times G \ar@<-.05cm>[r] \ar@<.05cm>[r] & X
}
\right\}
\simeq X/G\;,
\)
where we are now viewing the simplicial object on the left as a single smooth stack, with presheaves at each level. 

\begin{remark}
[Descent for quotient sheaves]
The equivalence \eqref{shdc} can be regarded as a descent condition for certain well-behaved quotient sheaves. In the proof of our main theorem, we will crucially use this identification to explicitly calculate a certain derived pullback which defines the inverse to the usual Pontrjagin-Thom construction. 
\end{remark}

There are several adjoint functors which relate the category ${\Sh}(\cartsp)$ to the 
category $\sset$ which will be useful for us throughout the paper. We now review these 
functors and the corresponding basic properties that we will need. 

\begin{remark} [Adjunctions for the model category of prestacks]
The $\infty$-category of smooth stacks participates in a 
quadruple adjunction $(\Pi\dashv \Gamma \dashv {\rm disc}\dashv {\rm codisc})$, in 
the set-up and notation of Schreiber \cite{Urs},
$$
\xymatrix{
\Sh(\cartsp)\ar@<.45cm>[rr]^-{\Pi}
\; \ar@<-.05cm>[rr]_-{\Gamma} && \; 
\infty-\mathscr{G}{\rm pd}
\ar@<-.05cm>[ll]_-{{\rm disc}}\ar@<.45cm>[ll]^-{\rm codisc}  \;,
}
$$
where 

\vspace{-3mm}
\begin{itemize}
\item $\Gamma$ is the \emph{global sections functor}, which evaluates a smooth stack 
on the point manifold.

\vspace{-3mm}
\item  ${\rm disc}$ is the \emph{`discrete' functor} which assigns each $\infty$-groupoid 
to the constant prestack and then stackifies the result.

\vspace{-3mm}
\item  $\Pi$ is the \emph{geometric realization functor} (which we 
will describe in further detail below). 
\end{itemize}

\vspace{-3mm}

\end{remark}

The functor in this adjunction of which we will make the most use is the geometric realization 
functor $\Pi$. Notice that it makes sense to take the locally constant stack associated to the 
simplicial set $\Pi(X)$. In fact, this operation is functorial and one gets a natural transformation
\(\label{shape modality}
\xymatrix{
\eta:{\rm id}\ar[r] & 
L({\rm const})\circ \Pi=:{\bf \Pi}
}\;,
\)
where $L$ denotes the \emph{stackification functor} already encountered above. 
The natural transformation ${\bf \Pi}$, sometimes called the \emph{shape} map,
can be regarded the map which sends a smooth stack to its associated homotopy 
invariant stack \cite{Urs}\cite{Lur}. 
This map will appear in our construction of the smooth cobordism category. It turns 
out that on smooth manifolds, this functor has a very convenient presentation. 
\begin{remark}[Singular realization]
\label{singular realization}
The following result, due to Pavlov \cite{Pa}, will be useful for us. 
Let ${\bf Sing}: {\Sh}(\mathscr{M}{\rm an})\longrightarrow {\PSh}(\mathscr{M}{\rm an})$ 
be the functor given by sending 
a smooth stack $F$, on the large site of smooth manifolds, to the
 homotopy colimit of smooth stacks (taken in the category of prestacks)
$$
{\bf Sing}(F):={\rm hocolim}_{[n]\in \Delta^{\rm op}}{\bf Map}(\Delta^n,F)\;,
$$
where $\Delta^n$ is the smooth $n$-simplex, viewed as a manifold with corners. 
Then ${\bf Sing}$ factors through ${\Sh}(\mathscr{M}{\rm an})$ and the canonical 
map ${\bf Sing}\to {\bf \Pi}$ is an equivalence \cite{Pa}.
\end{remark}

We use the above result to establish the following useful equivalences. 
This entails working with the category of compactly generated weak Hausdorff 
spaces (CGWH) \cite{Mc}. This category is well-behaved for mapping spaces
 and is closed under pushouts which involve inclusion of a closed subspace
 (see \cite{Str}).

\begin{proposition}
[Equivalences of realizations]
\label{eor pro}
For every smooth manifold $M$, and  smooth stack $X$ in $\sh_{\infty}(\cartsp)$, we have 
an equivalence of simplicial realizations

\vspace{-3mm}
$$
\Pi(\mathbf{Map}(M,X))\simeq \map(\Pi(M),\Pi(X))\;.
$$

\vspace{-2mm}
\noindent Moreover, we have an equivalence of topological realizations

\vspace{-3mm}
$$
\vert \Pi(\mathbf{Map}(M,X))\vert \simeq \map_{\rm CGWH}(M,\vert \Pi\vert (X))\;.
$$
\end{proposition}
\theproof
Since $\cartsp\subset \mathscr{M}{\rm an}$ is a dense subsite, the category $\sh_{\infty}(\cartsp)$ 
is equivalent to the category of smooth stacks on manifolds and this equivalence is exhibited by restriction to cartesian spaces. Thus, the claim will follow provided that 
it holds on the larger site of smooth manifolds. By Remark \ref{eor pro}, we can compute
$$
\Pi ({\bf Map}(M,X))\simeq \Gamma\circ {\bf Sing}({\bf Map}(M,X))\;.
$$
Calculating the homotopy colimit in simplicial sets via the diagonal of the relevant bisimplicial set, we see 
that the simplicial set on the right has simplices at level $k$ is given by the set of maps 
$M\times \Delta^k\to X_k$, where $X_k$ is the presheaf at level $k$ of $X$. The face maps induced 
by $\Delta^k$ and those of $X$. Using the closed Cartesian structure on ordinary presheaves, these 
maps are in bijective correspondence with the set of maps of the form $M\to {\bf hom}(\Delta_k,X_k)$,
and this correspondence clearly respects the face and degeneracy maps. This set is, in turn, immediately 
identified with the $k$-simplices of $\map(M,{\bf Sing}(X))$, taken in prestacks. Since ${\bf Sing}(X)$ 
satisfies descent, this gives an equivalence
$$
\Gamma\circ {\bf Sing}\big({\bf Map}(M,X)\big)\simeq \map(M,{\bf Sing}(X))\simeq \map(M,{\bf \Pi}(X))\;.
$$
But then, by the adjunction $\Pi\dashv {\rm disc}$, the space on the right is $\map(\Pi(M),\Pi(X))$. 
The second claim is well-known and follows from the Quillen equivalence between CGWH spaces 
and simplicial sets \cite{Qu}.
\endofproof

\begin{remark}
[Simplicial realization]
If $X$ is a smooth manifold, thought of as a smooth stack, remark \ref{eor pro} implies that we can identify the geometric realization as the simplicial set of $X$, which at level $n$ is given by the set $C^{\infty}(\Delta^n,X)$. Thus, we see that under the \emph{further} geometric realization $\vert \cdot \vert:\sset\to \mathscr{T}{\rm op}$, we recover the 
geometric realization of the smooth singular simplicial set. By the theorem of Milnor \cite{Mil}, this is indeed 
equivalent to $X$ itself.
\end{remark}

\medskip
The previous discussion gives the notion of geometric realization of a smooth stack, but there is also a notion of realization 
which is intrinsic to the category of smooth stacks. Indeed, we can form the \emph{stacky realization} of a simplicial diagram of stacks via the left Kan extension of the above inclusion along the Yoneda 
embedding $\Delta^{\rm op}\into [\Delta^{\rm op}, {\PSh}(\cartsp)]$. We can compute this Kan extension via the 
local formula which gives the \emph{classifying stack of a simplicial object}.
$$
\BB(X_{\bullet}):=\int^{[n]\in \Delta}X_n\times \underline{ \Delta} [n] \;,
$$
where $X_n$ is stack in the $n$-th level of $X_{\bullet}$, and the underline indicates that we are taking
the constant stack associated to the corresponding simplicial set. More explicitly, the coend can be computed as 
the coequalizer of the diagram
{\small
$$
\xymatrix{
\underset{d:[k]\to[m]}{\coprod} X_k\times \underline{ \Delta}[m]~
\ar@<.06cm>[rr]^-{{\rm id}\times d}\ar@<-.06cm>[rr]_-{d\times{\rm id}} 
&& ~\underset{[n]}{\coprod}X_n\times \underline{ \Delta}[n]\;,
}
$$
}
where $d:[k]\to[m]$ is any map in $\Delta$. In general, the result of this construction is not a smooth stack. Thus, we need to stackify the result so that we stay in the right category. In the local model structure on prestacks, these two objects are equivalent and we will often denote both the the prestack $\BB(X_{\bullet})$ and its stackification by the same symbol. 

\begin{proposition}
[Stacky realization as homotopy colimit]
\label{srhl}
The stacky realization $\BB(X_{\bullet})$ computes the homotopy colimit over $X_{\bullet}$ in $\sh_{\infty}(\cartsp)$.
\end{proposition}
\theproof
First observe that, before stackifying, the prestack $\BB(X_{\bullet})$ is objectwise the realization of the resulting bisimplicial set. Since the realization computes the homotopy colimit in $\sset$ and homotopy colimits in prestacks are computed objectwise, we see that $\BB(X_{\bullet})$ is a model for the homotopy colimit in prestacks. Homotopy colimits of diagrams of stacks, taken in prestacks, are always computed by stackifying the result. The claim follows immediately.
\endofproof

\begin{remark} 
[Stackifying homotopy colimits]
\label{pre to stack} 
Proposition \ref{srhl} invokes a trick that will be frequently used throughout this paper. 
Namely, that if we are given a diagram of stacks $F:J\to \Sh(\cartsp)$, the homotopy colimit can be 
computed by first computing the homotopy colimit in prestacks, and then the stackifying the result. This immediately follows from the fact that the stackification functor $L$ is left adjoint to the inclusion $i:\Sh(\cartsp)\into \PSh(\cartsp)$.
\end{remark}
\medskip

The following example will be of particular interest.

\begin{example}[Classifying stack of a smooth category]
Let $\mathscr{C}$ be an internal category of smooth stacks (i.e. a smooth category). This is, by definition, 
a diagram of the form
$
\xymatrix{
X   \ar@<-.07cm>[r] \ar@<.07cm>[r] &  Y
}
$
 along with a map
$$
m:X\times_{Y}X\longrightarrow X\;,
$$
which gives a composition law, and which makes the relevant diagrams commute. Taking the  
{\v C}ech nerve of this diagram gives a simplicial diagram of the form
$$
\xymatrix{
X\times_{Y}X\times_{Y} X
\ar@<.18cm>[r] \ar@<-.18cm>[r] \ar@<.065cm>[r]\ar@<-.065cm>[r]  
& X\times_{Y}X  
\ar@<.095cm>[r]^-{m} \ar[r] \ar@<-.095cm>[r]  & 
X   \ar@<-.06cm>[r] \ar@<.06cm>[r] &  Y\;.
}
$$
We will denote the stacky realization of the nerve by $\BB\mathscr{C}$ and call it the 
\emph{classifying stack of the smooth category} $\mathscr{C}$.
\end{example}

We now show that the geometric realization of the classifying stack of a simplicial object coincides with 
the usual notion of classifying space of the underlying topological space (viewed simplicially). 
Thus the geometric realization of the stacky realization of a smooth category 
is the topological realization of geometric realization of the category.
Note that the notion of a classifying space of a topological category is described in 
\cite{Segal}\cite{Mo}\cite{We}. 

\begin{proposition}
[Commuting realizations]
\label{Prop real}
Let $X_{\bullet}$ be a simplicial object in ${\Sh}(\cartsp)$. We have a homotopy equivalence of
CW-complexes
$$
\vert \BB(X_{\bullet}) \vert \cong {\rm B}(\vert X_{\bullet} \vert)\;,
$$
where $\vert \cdot \vert: {\Sh}(\cartsp)\to \sset\to \mathscr{T}{\rm op}$ denotes the composite of geometric realization 
functors and $\vert X_{\bullet}\vert$ is the simplicial object in $\mathscr{T}{\rm op}$ given by level-wise geometrically realizing. 
\end{proposition}
\theproof
First note that $\Pi$ and $\vert \cdot \vert$ both preserve (even strict) colimits and products. \footnote{The functor $\Pi$ is presented by the colimit operation on prestacks, which is left adjoint to the constant prestack functor.}
 Thus, we have a homeomorphism
{\small 
$$
\vert \BB(X_{\bullet}) \vert \cong  {\rm coeq} \Big\{\xymatrix{
\underset{d:[k]\to[m]}{\coprod} \vert X_k \vert \times \vert \underline{\Delta}[m]\vert  
\ar@<.06cm>[rr]^-{{\rm id}\times d}\ar@<-.06cm>[rr]_-{d\times{\rm id}} 
&& \underset{[n]}{\coprod} \vert X_n\vert \times \vert\underline{\Delta}[n]\vert
}\Big\}\;.
$$
}
Since $\underline{\Delta}[k]$ is locally constant, we have an equality of simplicial sets $\Pi(\underline{\Delta}[k])=\Delta[n]$, where $\Delta[n]$ is the simplicial $n$-simplex. The geometric realization of $\Delta[n]$ 
is equal to the standard topological $n$-simplex $\Delta^n$ (as follows from \cite[Theorem 1]{Mil}). Thus, 
the right hand side of the above homeomorphism is exactly the topological realization of 
the complex $\vert X_{\bullet} \vert$. Finally, since the stackification morphism
$L:\BB(X_{\bullet})\to  L(\BB(X_{\bullet}))$ is a weak equivalence in the local model structure on 
prestacks, and geometric realization preserves these equivalences, we have a weak homotopy equivalence
$$
\vert L \vert:\vert \BB(X_{\bullet}) \vert \longrightarrow
 \vert L(\BB(X_{\bullet})) \vert\;.
$$
Since $\vert \cdot \vert$ takes values in CW-complexes, the claim follows from Whitehead's Theorem.
\endofproof

\subsection{Classification of vector bundles in smooth stacks}
\label{Sec class}

In this section we will define a stacky operation on vector bundles which is analogous to the Thom space construction in topological spaces. We begin by describing a proper notion of a vector bundle in the setting of smooth stacks (see \cite{NSS}). This just places a familiar concept in a new, more general, setting.
\begin{definition}
[Vector bundle over a stack]
Let $X$ be a smooth stack. A \emph{vector bundle} with fiber $V$ is a fibration $\pi:\eta\to X$ 
between two smooth stacks such that for there is an effective epimorphism $p:Y\to X$, so that we have a homotopy Cartesian square
$$
\xymatrix@=1.5em{
V \times Y\ar[rr]\ar[d]_-{\rm proj} && \eta\ar[d]^-{\pi}
\\
Y\ar[rr]_-{p} && X\;.
}
$$ 
\label{Def vb}
\end{definition}
Definition \ref{Def vb} recovers the usual definition of a locally trivial vector bundle over a smooth manifold. Indeed, if $X$ is a smooth manifold then any good open cover $\{U_{\alpha}\}$ of $X$ defines an effective epimorphism $\coprod_{\alpha}U_{\alpha}\to X$. Taking the full nerve of the cover, we are led to the diagram
{\small
\(\label{iterated pullback}
\xymatrix@R=1.5em{
\hdots \ar@<-.1cm>[r]\ar[r]\ar@<.1cm>[r]  &
\coprod_{\alpha\beta} U_{\alpha\beta} \times V\ar@<-.05cm>[rr]\ar@<.05cm>[rr] \ar[d] 
&& \coprod_{\alpha} U_{\alpha} \times V \ar[rr]\ar[d] && \eta\ar[d]^-{\pi}
\\
\hdots \ar@<-.1cm>[r]\ar[r]\ar@<.1cm>[r] &
\coprod_{\alpha\beta} U_{\alpha\beta} \ar@<-.05cm>[rr]\ar@<.05cm>[rr] &&
 \coprod_{\alpha} U_{\alpha} \ar[rr] && X\;,
}
\)
}
where the first square is homotopy Cartesian. By the Pasting Lemma, the first square also implies that
 the homotopy fiber of $\eta\to X$ is isomorphic to $V$. 

\medskip
The above diagram may seem redundant, since the first pullback square seems to give all the information about
 the local triviality of the bundle. Note, however, that the higher stages in the simplicial diagram in the top tell 
 us how to glue together the total space of the bundle. The maps in this simplicial diagram are not just determined by the combinatorial data needed to glue together $X$, but also the transition functions of the bundle. The axiom of descent for an $\infty$-topos (see Rezk \cite{Re} and Lurie \cite{Lur}) implies that whenever we have an iterated pullback diagram, such as \eqref{iterated pullback}, then the first square is homotopy Cartesian if and only if the top right arrow is a homotopy colimiting cocone. Thus, we must have
{\small
$$
\xymatrix{
 {\rm hocolim}\Big\{ \hdots \ar@<-.1cm>[r]\ar[r]\ar@<.1cm>[r]  & \coprod_{\alpha\beta}U_{\alpha\beta} \times V\ar@<-.05cm>[r]\ar@<.05cm>[r] & \coprod_{\alpha}U_{\alpha}\times V \Big\} \ar[r]^-{\simeq} & \eta\;.
}
$$
}
which allows us to recover $\eta$ from local data. In the case where $\eta$ is zero-truncated, like a smooth manifold, this reduces to the usual construction for the total space of the bundle
{\small
$$
\xymatrix{
 {\rm colim} \Big\{ \coprod_{\alpha\beta} U_{\alpha\beta} \times V\ar@<-.05cm>[r]\ar@<.05cm>[r] & \coprod_{\alpha} U_{\alpha} \times V \Big\} \ar[r]^-{\cong} & \eta\;,
}
$$
}
\begin{remark} [General fibers and morphisms]
\noindent {\bf (i)} Definition \ref{Def vb} works just as well for bundles with any fibers, not just vector spaces. Thus, one could consider 
the more general theory of $\infty$-bundles, where the fibers are allowed to have higher simplicial data in smooth stacks.
For us, however, we will only be concerned with the zero-truncated fibers as our main interest
 is structures on the tangent bundle of a smooth manifold.

\noindent {\bf (ii)} One can also define a morphism of vector bundles $\eta\to \xi$ over $M$ in the 
obvious way: namely a morphism of smooth stacks which induces a linear isomorphism on the fibers.
 Since linear isomorphisms are in particular smooth, this condition is well-defined in smooth stacks.
\end{remark}

Let us first recall that there is a smooth stack which classifies locally trivial principal 
${\rm O}(d)$-bundles, denoted $\BB{\rm O}(d)$. Up to equivalence, this stack can be 
defined simply as the homotopy orbit stack of the smooth sheaf ${\rm O}(d)$ acting on
 the point manifold $\ast$.  However, we will need an explicit model of this stack with which 
 to work  in practice (see \cite{FSSt}\cite{SSS3}). 
 
\begin{definition}
[Classifying stack of orthogonal bundles]
We define the classifying stack of smooth, locally trivial, principal ${\rm O}(d)$-bundles as the
 stack given level-wise by
\( \label{orthogonal orbit stack}
\BB{\rm O}(d):=\left\{\xymatrix{
\hdots \ar@<.15cm>[r] \ar@<-.15cm>[r] \ar@<.05cm>[r]\ar@<-.05cm>[r]  &
 {\rm O}(d)\times {\rm O}(d) \ar@<.1cm>[r] \ar[r] \ar@<-.1cm>[r] & 
 {\rm O}(d) \ar@<-.05cm>[r] \ar@<.05cm>[r] & \ast
}\right\}\;,
\)
where the face maps $d_k:{\rm O}(d)^n\to {\rm O}(d)^{n-1}$ are given by the projections for $k\leq n$ and 
$d_{n+1}:{\rm O}(d)^n\to {\rm O}(d)^{n-1}$ sends $(g_1,g_2\hdots, g_n)\mapsto (1,\hdots, 1, g_1g_2\hdots g_n)$, for
$g_i\in C^{\infty}(-,{\rm O}(d))$. 
\end{definition}

For any Lie group $G$, one can define $\BB G$ analogously. Note that, since we are working over the small site of 
Cartesian spaces, the simplicial presheaf defined by \eqref{orthogonal orbit stack} \emph{does} in fact satisfy 
descent (see \cite{FSSt} for details) and is, therefore,
 a well-defined object in ${\Sh}(\cartsp)$. 

\medskip
This stack indeed classifies locally trivial, smooth principal ${\rm O}(d)$-bundles. In fact, in \cite{FSSt} it was shown 
that at the level of connected components we have a natural bijection
$$
\pi_0\map\big(\check{C}(\{U_{\alpha}\}),\BB{\rm O}(d)\big)\cong \check{H}^1(\{U_{\alpha}\},{\rm O}(d))\;,
$$
where on the left we have replaced the smooth manifold $M$ by the homotopy colimit over the {\v C}ech nerve of some good open cover 
\footnote{From the point of view of model category theory, these replacements are cofibrant replacements in the local projective model structure \cite{DI}.} and on the right we have the 
group of nonabelian {\v C}ech cohomology with respect to the cover. Such a class forms precisely the data 
needed to construct an isomorphism class of principal ${\rm O}(d)$-bundles, and each representative of the class corresponds to a choice of transition functions of the bundle. Even more generally, it was shown in \cite{FSSt}\cite{SSS3}\cite{Urs} that this isomorphism lifts to an equivalence of $\infty$-groupoids
\(
\map\big(M,\BB{\rm O}(d)\big)\overset{\simeq}{\longrightarrow}{\rm Bun}(M)\;,
\)
where on the right we have the $\infty$-groupoid of locally trivial principal ${\rm O}(d)$-bundles and on the left we have the derived mapping space between $M$ and $\BB{\rm O}(d)$.

\medskip
In the proof of the main theorem, we will need to have explicit cocycle representations for maps out of a smooth manifold to certain classifying stacks. We now illustrate how to derive the {\v C}ech cocycle data in this case. Let us take the model for $\check{C}(\{U_{\alpha}\})$ given by the realization over the nerve (see \cite{FSSt} for more on this argument). The resulting stack takes the form
{\small
$$
\xymatrix{
\hdots \ar@<.15cm>[r] \ar@<-.15cm>[r] \ar@<.05cm>[r]\ar@<-.05cm>[r]  & \;
\coprod_{\alpha\beta\gamma }U_{\alpha\beta\gamma}\times \underline{\Delta}[2]
\; \ar@<.1cm>[r] \ar[r] \ar@<-.1cm>[r] & \;
\coprod_{\alpha\beta}U_{\alpha\beta}\times \underline{\Delta}[1] 
\;\ar@<-.05cm>[r] \ar@<.05cm>[r] & \;
 \coprod_{\alpha}U_{\alpha}\times \underline{\Delta}[0]\;,
}
$$
}
where the face maps are induced by the inclusion and the usual face maps of $\underline{\Delta}[k]$. Since $\BB{\rm O}(d)$ is 1-truncated, the Yoneda Lemma along with descent implies that a map from this stack to $\BB{\rm O}(d)$ is equivalently given by a commutative diagram
{\small
$$
\xymatrix@R=1.3em{
\underline{\Delta}[2]\ar[rr] && \prod_{\alpha\beta\gamma }\BB{\rm O}(d)(U_{\alpha\beta\gamma})
\\
\underline{\Delta}[1]\ar@<-.1cm>[u]\ar[u]\ar@<.1cm>[u]\ar[rr] && 
\prod_{\alpha\beta }\BB{\rm O}(d)(U_{\alpha\beta})\ar@<-.1cm>[u]\ar[u]\ar@<.1cm>[u]
\\
\underline{\Delta}[1]\ar@<-.05cm>[u]\ar@<.05cm>[u]\ar[rr] && 
\prod_{\alpha }\BB{\rm O}(d)(U_{\alpha})\ar@<-.05cm>[u]\ar@<.05cm>[u] \;,
}
$$
}
where the vertical maps on the right are given by restriction and those on the left are given by the face inclusions. From the definition of $\BB{\rm O}(d)$, one sees that such a diagram determines -- and is uniquely determined
-- by a choice of smooth ${\rm O}(d)$-valued function $g_{\alpha\beta}:U_{\alpha\beta}\to {\rm O}(d)$ on intersections satisfying the cocycle condition
$
g_{\alpha\beta}g_{\beta\gamma}g_{\gamma\delta}=1
$
on 3-fold intersections $U_{\alpha\beta\gamma}$. It follows immediately that we have an isomorphism
$$
\hom\big(\check{C}(\{U_{\alpha}\}),\BB{\rm O}(d)\big)
\cong 
\check{C}^{1}\big(\{U_{\alpha}\},{\rm O}(d)\big)\;,
$$
where the right hand side is the set of degree 1 nonabelian {\v C}ech cocycles with values in ${\rm O}(d)$. 

\medskip
In the same way that locally trivial ${\rm O}(d)$-bundles are classified by maps to $\BB{\rm O}(d)$, there is a universal bundle $\mathbfcal{U}(d)\to \BB{\rm O}(d)$ which classifies vector bundles over $M$, via pullback along a classifying map. Indeed, 
let $\mathbfcal{U}(d)$ be the smooth stack defined by the homotopy quotient $\RR^d/\!/{\rm O}(d)$, where ${\rm O}(d)$ acts in the usual way. Lemma \ref{orbit stack} provides us with an  identification
\( \label{universal bundle}
\mathbfcal{U}(d)\simeq \left\{\xymatrix{
\hdots \ar@<.15cm>[r] \ar@<-.15cm>[r] \ar@<.05cm>[r]\ar@<-.05cm>[r]  &
 \RR^d\times {\rm O}(d)\times {\rm O}(d) \ar@<.1cm>[r] \ar[r] \ar@<-.1cm>[r] & 
 \RR^d\times {\rm O}(d) \ar@<-.05cm>[r] \ar@<.05cm>[r] & \RR^d
}\right\}\;,
\)
where the face maps $d_k:{\rm O}(d)^n\times \RR^d \to {\rm O}(d)^{n-1}\times \RR^d$ are given by the projections for $k\leq n$ and 
$d_{n+1}:{\rm O}(d)^n\times \RR^d\to {\rm O}(d)^{n-1}\times \RR^d$ sends $(g_1,g_2\hdots, g_n,v)\mapsto (1,\hdots, 1, g_1g_2\hdots g_nv)$. 
 If we are given a map $f:M\to \BB{\rm O}(d)$ classifying an ${\rm O}(d)$-bundle on $M$, 
 then the pullback
\(\label{pullback bundle}
\xymatrix@R=1.5em{
\eta \ar[rr]^{g} \ar[d]_{\pi}&& \mathbfcal{U}(d) \ar[d]^-{\pi}
\\
M \ar[rr]^-{f} && \BB{\rm O}(d)  
}
\)
is a locally trivial vector bundle over $M$. The trivialization comes from the fact that when we replace $M$ by a {\v C}ech resolution $\check{C}(\{U_{\alpha}\})$. One such pullback is explicitly computed as $\check{C}(\{U_{\alpha}\times \RR^d\})$, where the face maps in the {\v C}ech nerve diagram depend on the transition functions determined by $f$. Since any two pullbacks are canonically isomorphic, this resolution defines a local trivialization of the bundle. 

\medskip
Notice that immediately one obtains that $\mathbfcal{U}(d)$ defines a vector bundle in the sense of

\begin{remark} [Choice of {\v C}ech resolution]
Note that it is not quite right to say that a map $M\to \BB{\rm O}(d)$ determines a locally trivial principal ${\rm O}(d)$-bundle on $M$. Strictly speaking, we need to first choose a {\v C}ech resolution of $M$ and then take a map from the resolution to $\BB{\rm O}(d)$ to get such a bundle. In general, the local structure depends on the choice of resolution, but the full mapping spaces out of different resolutions will be equivalent (this follows from general properties of simplicial model categories and the characterization of cofibrant objects via \cite{DI}). Henceforth, whenever we write a map out of a smooth manifold $M$, we will implicitly assume a choice of {\v C}ech resolution. If a specific choice is needed, we will make it clear.
\end{remark}

The picture is even more clear when one looks at the {\v C}ech cocycles obtained by mapping into $\mathbfcal{U}(d)$. Computing the cocycles in the same way as before, we see that for a good open cover $\{U_{\alpha}\}$ of $M$, a map to $\mathbfcal{U}(d)$ determines and is uniquely determined by a commutative diagram
 {\small 
 \(\label{cech data for universal}
\xymatrix@R=1.3em{
\underline{\Delta}[2]\ar[rr] && \prod_{\alpha\beta\gamma }
\mathbfcal{U}(d)(U_{\alpha\beta\gamma})
\\
\underline{\Delta}[1]\ar@<-.1cm>[u]\ar[u]\ar@<.1cm>[u]\ar[rr] &&
 \prod_{\alpha\beta }\mathbfcal{U}(d)(U_{\alpha\beta})\ar@<-.1cm>[u]\ar[u]\ar@<.1cm>[u]
\\
\underline{\Delta}[0]\ar@<-.05cm>[u]\ar@<.05cm>[u]\ar[rr] &&
 \prod_{\alpha }\mathbfcal{U}(d)(U_{\alpha})\ar@<-.05cm>[u]\ar@<.05cm>[u]\;.
}
\)
}
 By the definition of $\mathbfcal{U}(d)$, such a diagram is equivalently a choice of local section $s_{\alpha}:U_{\alpha}\to \RR^d$ and transition functions $g_{\alpha\beta}:U_{\alpha\beta}\to {\rm O}(d)$ such that 
 $s_{\alpha}=g_{\alpha\beta}s_{\beta}$
 on intersections and 
 $g_{\alpha\beta}g_{\beta\gamma}g_{\gamma\delta}=1$
 on 3-fold intersections. This is precisely the data needed to construct a global section of the bundle with transition functions $g_{\alpha\beta}$. At an abstract level, this is clear since the universal property of the pullback implies that if we are given a map $f:M\to \mathbfcal{U}(d)$ which fills the diagram \eqref{pullback bundle}, then we have an induced map $s:M\to \eta$ such that $\pi s={\rm id}$, i.e., a global section.

\subsection{The Grassmanian and its canonical bundles in the stacky setting}
\label{Sec gras}

As in the Madsen-Tillman construction \cite{MT}, we now turn our attention to embedded submanifolds and their tangent and normal bundles. Let us recall the Grassmannian manifold, whose underlying set is defined as the collection of $d$-planes in $\RR^{d+N}$. This set is given smooth structure via the quotient 
\(\label{gr man}
{\rm Gr}(d,N)={\rm O}(d+N)/{\rm O}(d)\times {\rm O}(N)\;,
\)
where ${\rm O}(d)\times {\rm O}(N)$ acts via matrix multiplication by block matrices. Similarly, we recall the definition of the real Steifel manifold ${\rm V}(d,N)$ as the set of orthonormal $d$-frames in $\RR^N$:
\(\label{st man}
{\rm V}(d,N)={\rm O}(d+N)/{\rm O}(d)\;.
\)
There are two canonical bundles over the Grassmannian: the \emph{tautological bundle}, which maps the points on a $d$-plane to the corresponding $d$-plane in $\RR^{d+N}$, and the \emph{orthogonal complement bundle}, given by mapping the points on the complement of a $d$-plane to the underlying $d$-plane.  Both bundles admit the structure of smooth manifolds and are identified, respectively, by the quotients
\(\label{uo man}
\mathcal{U}(d,N):=\RR^d\times_{{\rm O}(d)}{\rm V}(d,N)
\qquad \text{and} \qquad
\mathcal{U}^{\perp}(d,N):=\RR^N\times_{{\rm O}(N)}{\rm V}(N,d)\;.
\)
Here we think of $d$ as indexing the tangent direction and $N$ as indexing the normal direction, although clearly we could reverse the roles of $d$ and $N$, which would lead to the canonical  identification $\mathcal{U}(d,N)\cong \mathcal{U}^{\perp}(N,d)$. In each example, we can take the colimit over the canonical inclusion maps as $N\to \infty$. This leads to the corresponding infinite-dimensional manifolds ${\rm Gr}(d,\infty)$, ${\rm V}(d)$, $\mathcal{U}(d)$ and $\mathcal{U}^{\perp}(d)$. 

\medskip
The purpose of this section is to provide convenient descent data for each of these smooth manifolds, viewed as living in the category of smooth stacks. Notice that each of these smooth manifolds arise from quotients of Lie group actions and we can, therefore, use the identification \eqref{shdc} for each of these objects. Since we will use this identification frequently, we will spell it out in detail.

\begin{lemma}
[Homotopy orbit stack of a Lie group action]
\label{orbit stack}
Let $G$ be a Lie group, regarded as a zero-truncated smooth stack (i.e. a smooth sheaf). Suppose $G$ 
acts on a freely on a sheaf $X$ and suppose moreover that the corresponding quotient presheaf $X/G$ satisfies descent (i.e. is a sheaf). Then the prestack
\(\label{explicit classifying stack}
Y:=\left\{\xymatrix{
\hdots \ar@<.15cm>[r] \ar@<-.15cm>[r] \ar@<.05cm>[r]\ar@<-.05cm>[r]  
& X\times G\times G \ar@<.1cm>[r] \ar[r] \ar@<-.1cm>[r] 
& X\times G \ar@<-.05cm>[r] \ar@<.05cm>[r] & X
}\right\}\;,
\)
satisfies descent and we have an equivalence of smooth stacks
$$q:Y\overset{\simeq}{\longrightarrow} X/G\;,$$
where $q$ is the strict quotient map.
\end{lemma}
\theproof
Homotopy colimits of prestacks can be calculated objectwise (this follows, for example, from \cite[Theorem 4.2.4.1]{Lur}). Therefore, the homotopy colimit can be calculated as the presheaf which assigns to each $U\in \cartsp$, 
the diagonal of the bisimplicial set given by considering \eqref{explicit classifying stack} as a simplicial object in prestacks and evaluating at $U$. Since this diagram is constant in one simplicial direction, the prestack given by \eqref{explicit classifying stack} is a model for the homotopy orbits, taken in prestacks. Now we use the trick discussed in Remark \ref{pre to stack} to conclude that the stackification models the homotopy quotient in stacks. 

On the other hand, since $G$ acts freely on $X$, the stabilizer subgroups of $G$ must be trivial. But these subgroups are precisely the presheaves $\pi_1(Y,g)$, where $g:\ast\to G$ is a section. Since the higher simplicial homotopy groups vanish ($G$ is zero-truncated), we see that the map
$
q:Y\to X/G
$,
induced by the unit of the adjunction $\pi_0\dashv sk_0$, induces an isomorphism on all presheaves of homotopy groups. Since $X/G$ already satisfies descent and is equivalent to $Y$ at the level of prestacks, we immediately deduce that $Y$ satisfies descent and that $q$ is an equivalence of smooth stacks.
\endofproof

\begin{remark}
In what follows we will be working with various models for certain smooth stacks. In particular, we will occasionally want to distinguish between strict sheaves and smooth stacks which are equivalent to these strict sheaves. Whenever we are working only up to equivalence, or choosing models for smooth sheaves which have higher simplicial data, we will denote the corresponding stack by bold characters (for example, ${\bf Gr}(d,N)$). When we are identifiying the stack with a strict smooth sheaf, we will use upright roman characters (for example, ${\rm Gr}(d,N)$). This convention will be used throughout the remainder of the paper.
\end{remark}

Notice that in each of the above examples in expressions \eqref{gr man}, \eqref{st man} and \eqref{uo man}, the given 
action is free. Thus, Lemma \ref{orbit stack} applies and gives us the convenient descent data we have been looking for in each example. This now motivates us to define various stacks which will replace their classical counterparts. Lemma \ref{orbit stack} will imply that each of these objects is equivalent to the strict quotient and serves as another model for the \emph{homotopy} quotient in stacks. 
\begin{definition}
[Smooth Grassmannian stack]
\label{Def Gra}
We define the \emph{smooth Grassmannian stack} to be the smooth stack given level-wise by
$$
\hspace{-.1cm}
\mathbf{Gr}(d,N):=\Big\{
\xymatrix@=1.45em{
\hdots \ar@<.15cm>[r] \ar@<-.15cm>[r] \ar@<.05cm>[r]\ar@<-.05cm>[r] &  \ar@<.1cm>[r] \ar[r] \ar@<-.1cm>[r] ({\rm O}(d)\times {\rm O}(N))^2 \times {\rm O}(d+N)  & ({\rm O}(d)\times {\rm O}(N))\times {\rm O}(d+N) \ar@<-.05cm>[r] \ar@<.05cm>[r] & {\rm O}(d+N)
}\Big\}
$$
where the face maps $d_k$ are projections for $k\leq n$ and $d_{n+1}$ is given by the action of $({\rm O}(d)\times {\rm O}(N))$ 
on ${\rm O}(d+N)$ via the canonical inclusion ${\rm O}(d)\times {\rm O}(N) \into {\rm O}(d+N)$ as block matrices for $k=n+1$. 
\end{definition}

 
  We can also define the stacky analogue of the real Stiefel manifold.
 
\begin{definition}
 [Smooth Stiefel stack] Define the smooth Stiefel stack via
$$
\mathbf{V}(d,N):=\left\{\xymatrix@=1.45em{
\hdots \ar@<.15cm>[r] \ar@<-.15cm>[r] \ar@<.05cm>[r]\ar@<-.05cm>[r] &  \ar@<.1cm>[r] \ar[r] \ar@<-.1cm>[r] {\rm O}(N)^2 \times {\rm O}(d+N)  & {\rm O}(N)\times {\rm O}(d+N) \ar@<-.05cm>[r] \ar@<.05cm>[r] & {\rm O}(d+N)
}\right\} \;,
$$
where for $k\leq n$, the face maps $d_k$ are again projections and where $d_{n+1}$ is given by the action by matrix multiplication 
after the inclusion\footnote{We will use the notation $i^p_n$ to denote the embedding of $p \times p$ matrices as the upper left block in 
$n \times n$ matrices ($n \geq p$), and similarly $i_{q,n}$ for the embedding of $q \times q$ matrices as the 
lower right block in $n \times n$ matrices ($n \geq q$).} 
$i^{N}_{d+N}:{\rm O}(N)\into {\rm O}(d+N)$, which sends an orthogonal 
$(N\times N)$-matrix $Q$ to the block matrix with $Q$ in the upper left corner, the $d \times d$ identity matrix $1_{d}$ in the lower right corner and zeros everywhere else. 
\end{definition}

As in the case of the Grassmannian, this prestack satisfies descent, is zero-truncated and is equivalent to its 
classical counterpart ${\rm V}(d,N)$. Notice also that the smooth Stiefel stack admits an action of 
${\rm O}(d)$, given by the inclusion $i_{d,d+N}:{\rm O}(d)\into {\rm O}(d+N)$, which includes as a block matrix in the \emph{lower right corner}.

\begin{remark} 
[Non-contractibility of the Stiefel Stack in the limit]
Notice that $\mathbf{V}(d,N)$ is clearly \emph{not} contractible as a smooth stack as $N\to \infty$. Indeed, if it were then $\colim_{N}{\rm V}(d,N)$ would have to be \emph{diffeomorphic} to the point manifold  $\ast$ as $N\to \infty$. \footnote{Recall that $V(d,N)$ is a zero-truncated smooth manifold.} This is obviously not the case. For example, we can choose a sequence $Q_{N}$ of orthogonal matrices of the form
$$
Q_N=\underbrace{
\left(
\begin{array}{c|c}
1 & \ast \\
\hline
\ast & Q
\end{array}
\right)
}_{N+d}\;,
$$
%
with $1\neq Q\in {\rm O}(d)$. Clearly each $Q_N$ is a non identity element in ${\rm O}(N+d)/{\rm O}(N)$, and $Q_{N}\mapsto Q_{N+1}$, under the induced map ${\rm V}(d,N)\to {\rm V}(d,N+1)$. Hence, in the colimit $N\to \infty$, $Q_N$ converges to a nonidentity element $Q_{\infty}\in {\rm O}/{\rm O}$, since for no value of $N$ is $Q_N=1$. For this reason, we see that the stacks ${\bf Gr}(d,N)$ cannot approximate $\BB{\rm O}(d)$ as $N\to \infty$, \footnote{In fact, we already knew this, since ${\bf Gr}(d,N)$ is zero-truncated.} which is of course different than the case in topological spaces.
\end{remark}

\medskip
We can now define the universal orthogonal complement bundle 
over the Grassmannian. Indeed, through the map $i^{N}_{d+N}:{\rm O}(N)\into {\rm O}(d+N)$, 
the group ${\rm O}(N)$ acts on the orthogonal complement of $\RR^d$ in $\RR^{d+N}$ and this 
leads to the following definition.
\begin{definition}
[Universal vector bundle and it complement over Grassmannian stacks]
 \label{universal bundle}
{\bf (i)} We define the \emph{universal vector bundle} $\pi:\mathbfcal{U}(d,N)\to 
{\bf Gr}(d,N)$ over the smooth Grassmannian, by setting
$$
\mathbfcal{U}(d,N):=\big(\RR^d\times {\rm V}(d,N)\big)/\!/{\rm O}(d)\;,
$$
where the action of ${\rm O}(d)$ on $\RR^d$ is given by the usual action by linear maps, and on the 
product diagonally (with the action on the second factor as above). 

\item {\bf (ii)}
Similarly, we define the \emph{universal orthogonal complement bundle} $\pi^{\perp}:\mathbfcal{U}^{\perp}(d,N)\to {\bf Gr}(d,N)$ over the Grassmannian, by setting 
$$
\mathbfcal{U}^{\perp}(d,N):=\big((\RR^d)^{\perp}\times {\rm V}(N,d)\big)/\!/{\rm O}(N)\;,
$$
where the orthogonal complement $(\RR^d)^{\perp}$ of $\RR^d$ is taken in $\RR^{d+N}$ (and can, therefore, be identified 
with $\RR^N$) and the action of ${\rm O}(N)$ is given by the diagonal action. The bundle map is induced by projecting out the first factor
$(\RR^d)^{\perp}$.
\end{definition}

Since each of the above stacks is equivalent to its classical analogue, we immediately deduce the following.

\begin{proposition}
[Classical correspondence]
The geometric realization of all the above stacks are weak homotopy equivalent to their classical counterparts.
\end{proposition}\label{realization of stacks}


As in the case of the usual Grassmannian, we have obvious maps 
$\mathbf{Gr}(d,N)\into \mathbf{Gr}(d,N+1)$ which are induced by the inclusion $i^{N}_{1+N}:{\rm O}(N)\into {\rm O}(N+1)$. In accord with the classical case, we have the following splitting.

\begin{proposition} 
[Splitting of the universal orthogonal complement bundle]
\label{pullback trivializes}
 The homotopy pullback of $\mathbfcal{U}^{\perp}(d,N+1)$ to ${\bf Gr}(d,N)$ can be identified with the sum $\mathbfcal{U}^{\perp}(d,N)\oplus {\bf 1}$.
\end{proposition}
\theproof
This amounts to a direct calculation of the pullback of the bundle. The homotopy pullback can be computed degreewise in prestacks, and we see that in degree $k$ we are led to the homotopy Cartesian square
{\small
$$ 
\xymatrix@R=1.5em{
(\RR^d)^{\perp}\oplus \RR \times {\rm O}(d+N)\times \big({\rm O}(d)\times {\rm O}(N)\big)^k 
\ar[rr]\ar[d]  && 
(\RR^d)^{\perp} \times {\rm O}(d+N+1)\times \big({\rm O}(d)\times {\rm O}(N+1)\big)^k \ar[d]
\\
 {\rm O}(d+N)\times \big({\rm O}(d)\times {\rm O}(N)\big)^k
 \ar[rr] && 
 {\rm O}(d+N+1)\times \big({\rm O}(d)\times {\rm O}(N+1)\big)^k\;.
}
$$
}
Note that the right vertical maps project out the Euclidean spaces, and we have used the fact that the map on the right is an objectwise Kan fibration and all objects are objectwise fibrant in order to compute the homotopy pullback as the strict pullback. \footnote{Note that the orthogonal complements are taken in different dimensions.} Since the action of ${\rm O}(N)$ does not act on the factor $\RR$ in the upper left corner, we see that the resulting prestack is simply $({\rm Gr}(d,N)\times \RR)\oplus \mathbfcal{U}^{\perp}(d,N)$. Finally, since the stackification functor is left exact, it preserves homotopy Cartesian squares and products. We compute the pullback in stacks via stackification, which leads to the result.
\endofproof

We will similarly need the colimit of the orthogonal complement bundle. 

\begin{definition} 
[Universal orthogonal complement bundle]
We define the \emph{universal orthogonal complement bundle} 
$\mathbfcal{U}^{\perp}(d)\to {\bf Gr}(d,\infty)$ of the classifying stack of orthogonal 
bundles as the universal map induced by passing to colimits in Definition \ref{universal bundle}, i.e. 
$\mathbfcal{U}^\perp (d):=\colim_{k}\mathbfcal{U}^\perp(d, N)$. 
\end{definition}

Note that one can obtain {\v C}ech cocycle data for the bundles $\mathbfcal{U}^{\perp}(d)$ and $\mathbfcal{U}^{\perp}(d,N)$ in a completely analogous way to that of $\mathbfcal{U}(d)$, discussed in Section \ref{Sec class}. 

\begin{remark}
[Stacks vs. spaces]
 We could have deduced many of the propositions from classical results. We have chosen to illustrate how to work with these smooth stacks, since we will be using similar techniques in the proof of the main theorem.
\end{remark}

\section{Thom stacks and their smooth motivic spectra}
\label{Sec Thom}

\subsection{Thom stacks} 
\label{Sec Th}

Having defined vector bundles in the context of stacks, now we consider the stacky analogue of the 
corresponding Thom spaces. We begin with a classical discussion. Suppose $V$ 
is finite-dimensional vector space and equip $V$ with an inner product. Let $\pi:\eta\to M$ be a vector bundle with fiber $V$ over a smooth manifold $M$. Then $\eta$ inherits a metric from $V$ via locally trivializing patches and the transition functions of the bundle can be chosen to be orthogonal, via the metric. Let $D(V)$ denote the closed unit disc in $V$ with respect to this inner product. Consider the diagram
{\small
$$
\xymatrix@=1.5em{
\coprod_{\alpha\beta}U_{\alpha\beta}\times D(V)\ar@<.05cm>[rr] \ar@<-.05cm>[rr]\ar[d] &&\coprod_{\alpha}U_{\alpha}\times D(V)\ar[rr]\ar[d] && D(\eta)\ar[d]
\\
\coprod_{\alpha\beta}U_{\alpha\beta}\times V \ar@<.05cm>[rr] \ar@<-.05cm>[rr]  && 
\coprod_{\alpha}U_{\alpha}\times V\ar[rr] && \eta
\;,
}
$$
}
where the top maps are induced by the transition functions of the bundle on intersections and $D(\eta)$ denotes the unit disc bundle. All the squares to the left of the right most square are Cartesian and by the axiom of descent for the category of sheaves, the top diagram is a coequalizer precisely when the first square is Cartesian. This ensures that we do have a 
well-defined disc bundle with trivializing charts which are compatible with the ones for $\eta$. We also have a similar picture for the sphere bundle $S(\eta)$ and we have a canonical map $S(\eta)\into D(\eta)$.

\medskip
This discussion works well for total spaces which are sheaves, but generally we will need to consider total spaces which are smooth stacks. The same construction works in this case too, but we need to extend the diagram on the left via the entire hypercover
{\small
\(\label{simplicial diagram disc}
\xymatrix@=1.5em{
\ar@<.1cm>[r]\ar[r] \ar@<-.1cm>[r]\hdots & \coprod_{\alpha\beta}U_{\alpha\beta}\times D(V)
\ar@<.05cm>[rr] \ar@<-.05cm>[rr]\ar[d] && 
\coprod_{\alpha}U_{\alpha}\times D(V)\ar[rr]\ar[d] && D(\eta)\ar[d]
\\
\ar@<.1cm>[r]\ar[r] \ar@<-.1cm>[r]\hdots & \coprod_{\alpha\beta}U_{\alpha\beta}\times V 
\ar@<.05cm>[rr] \ar@<-.05cm>[rr]  && \coprod_{\alpha}U_{\alpha}\times V\ar[rr] && \eta
\;.
}
\)
}
Then $\eta$ will be the \emph{homotopy} colimit over the bottom simplicial diagram and $D(\eta)$ will be the homotopy colimit over the top simplicial diagram. Again the axiom of descent implies that all squares are homotopy Cartesian and the top and bottom maps are homotopy quotient maps. This way, we have a well-defined notion of a unit disc and sphere bundle for vector bundles which have smooth stacks as total spaces. 

\medskip
In fact, we can continue the discussion and consider the situation in even greater generality. Indeed, we can perform the same construction for vector bundles with total spaces given by smooth stacks and base spaces given by smooth stacks. We simply replace our cover of the manifold $M$ by effective epimorphisms $Y\to X$. 

\begin{definition}
[Thom stack]
\label{DefThomstack}
Let $V$ be a finite-dimensional vector space, equipped with inner product. For any smooth stack $X$, 
we define the {\rm Thom stack} of a $V$-bundle, $\pi:\eta\to X$ to be the homotopy quotient stack 
$$
{\rm Th}(\eta):=D(\eta)/\!/S(\eta)\;.
$$
\end{definition}

If the map $S(\eta)\into D(\eta)$ is an objectwise cofibration (i.e. monomorphisms of simplicial sets), then the homotopy quotient in prestacks $D(\eta)/\!/S(\eta)$ can be modeled by the strict quotient. Thus, the stackification of the strict quotient is a model for the Thom stack in this case. The resulting stack is, in fact, a sheaf up to equivalence.
It is the diffeological space whose plots are smooth functions to the
$n$-disk bundle, where all plots which land entirely in the boundary of
the $n$-disk bundle get identified.

\begin{remark}
[`Stacky sphere']
If ${\rm dim}(V)=n$, the quotient stack $D(V)/S(V)$ will serve as a model for the $n$-dimensional sphere and we will denote this stack accordingly as $D^n/\partial D^n$. 
However, it is important to note that this \emph{is not} the stack represented by the sheaf of smooth plots of the unit $n$-sphere in $V\oplus \RR$, nor is it the simplicial sphere $S^n:=\Delta[n]/\partial\Delta[n]$.  
\end{remark}

We immediately have the following. 

\begin{proposition}
[Thom stack of a trivial bundle]
\label{sphere wedge}
Let ${\bf n}:\RR^n\times X\to X$ be the trivial rank $n$ bundle on a smooth stack $X$. We have an 
equivalence
$$
{\rm Th}({\bf n})\simeq D^n/\partial D^n\wedge X_+\;.
$$
\end{proposition}
\theproof
This is immediate from the definition of the Thom stack.
\endofproof

As in the classical case, we would like the Thom stacks to behave well with respect to the smash product of pointed stacks. With an appropriate understanding of the $n$-sphere, we see that this is indeed the case. 
\begin{proposition}
[Properties of Thom stacks]
\label{Pro pro}
Let $\xi \to X$ and $\eta \to Y$ be vector bundles of rank $n$ and $m$, respectively. The Thom stack operation ${\rm Th}$ 
 satisfies the following properties.

\noindent {\bf (i)} {\rm (Functorial)}  Given a morphism of bundles $\phi:\xi\to \eta$ which restricts to an isometric injective linear map on the fibers, there is an induced
 morphism of Thom stacks ${\rm Th}(\phi):{\rm Th}(\xi)\to {\rm Th}(\eta)$.

\noindent {\bf (ii)} {\rm (Multiplicative)} We have an equivalence of smooth stacks
$
{\rm Th}(\xi\square \eta)\simeq {\rm Th}(\xi)\wedge {\rm Th}(\eta)
$,
where $\xi\square \eta\to X\times Y$ is the external sum bundle.

\label{Prop external}
\end{proposition}
\theproof
Part {\bf (i)} follows directly from the commutative diagram
$$
\xymatrix@=1.5em{
S(\xi) \ar@{^{(}->}[d] \ar[rr]^-{\phi} && S(\eta) \ar@{^{(}->}[d]
\\
D(\xi)\ar[rr]^-{\phi} && D(\eta)
\;.
}
$$
To prove part {\bf (ii)},  note that we have a canonical map
\(\label{quotient}
{\rm Th}(\eta\square \xi)\longrightarrow {\rm Th}(\eta)\wedge {\rm Th}(\xi)
\)
induced by the composition $D(\xi\times \eta)\into D(\xi)\times D(\eta)\to 
{\rm Th}(\xi)\wedge {\rm Th}(\eta)$, which descends to the quotient. This map is an 
equivalence fiberwise since we have an isomorphism of pointed sheaves 
\footnote{This is essentially by inspection, along with the fact that the the interior of 
$D^{n+m}$ is diffeomorphic to the product of the interiors of $D^n$ and $D^m$.}
$$
D^{n+m}/\partial D^{n+m}\cong D^{n}/\partial D^{n}\wedge D^{m}/\partial D^{m}\;.
$$ 
In ${\Sh}(\cartsp)$, maps which are surjective of sheaves of connected components and 
which induce fiberwise equivalences are equivalences of smooth stacks (this follows, 
for example, from the long exact sequence on sheaves of homotopy groups). It follows 
that the map \eqref{quotient} is an equivalence.
\endofproof

The next result allows us to identify the Thom stack of the various universal bundles over the 
Grassmannian with particularly nice smooth stacks.
\begin{proposition}
[Thom stacks of universal bundles]
The total stacks of the Thom stacks of the universal bundles 
$\mathbfcal{U}(d)\longrightarrow \BB{\rm O}(d)$, 
$\mathbfcal{U}(d,N)\longrightarrow {\bf Gr}(d,N)$ and 
$\ \mathbfcal{U}^{\perp}(d,N)\longrightarrow{\bf Gr}(d,N)$ can be identified as

\vspace{-8mm}
\bea
{\rm Th}(\mathbfcal{U}(d)) &\simeq & (D^d/\partial D^d)/\!/{\rm O}(d)\;,
\\
{\rm Th}(\mathbfcal{U}(d,N))&\simeq  & \big(D^d/\partial D^d\wedge {\rm V}(d,N)_+\big)/\!/{\rm O}(d)\;,
\\
{\rm Th}(\mathbfcal{U}^{\perp}(d,N))&\simeq & \big(D^N/\partial D^N\wedge {\rm V}(N,d)_+\big)/\!/{\rm O}(N)\;,
\eea
where, in each example, the action on $D^d/\partial D^d$ is inherited from the action of ${\rm O}(d)$ on $\RR^d$. 
\footnote{Note that the unit disc and unit sphere are fixed under the action of ${\rm O}(d)$. Thus, we indeed 
get a well-defined action on the quotient.}
\end{proposition}
\theproof
We will prove the claim for the first Thom stack and the others are proved analogously. 
We first make the identification in prestacks and then stackify the result. In prestacks, the map
$$
\xymatrix{
\partial D^d/\!/{\rm O}(d) \; \ar@{^{(}->}[r] & D^d/\!/{\rm O}(d)
}
$$
is an objectwise monomorphisms. Thus, the homotopy quotient can be computed as the strict 
quotient of prestacks. This, in turn, is immediately identified with $(D^d/\partial D^d)/\!/{\rm O}(d)$. 
Using our standard trick for computing homotopy quotients in smooth stacks gives the identification.
\endofproof

Notice that for the smooth stacks ${\rm Th}(\mathbfcal{U}(d,N))$ and ${\rm Th}(\mathbfcal{U}^{\perp}(d,N))$, we again have equivalences 
$$
{\rm Th}(\mathbfcal{U}(d,N))\simeq {\rm Th}(\mathcal{U}(d,N))
\qquad \text{and} \qquad
 {\rm Th}(\mathbfcal{U}^{\perp}(d,N))\simeq {\rm Th}(\mathcal{U}^{\perp}(d,N))\;,
$$
where on the right we have the corresponding classical Thom stacks given by applying the usual Thom space construction to the smooth manifolds $\mathcal{U}(d,N)$ and $\mathcal{U}^{\perp}(d,N)$ in the ambient category of smooth stacks. All of these constructions geometrically realize to give the underlying Thom spaces.

\begin{proposition}
[Classical correspondence]
\label{geometric realization of universal bundle}
The geometric realization of the Thom stacks ${\rm Th}(\mathbfcal{U}(d))$, ${\rm Th}(\mathbfcal{U}(d,N))$ and ${\rm Th}(\mathbfcal{U}^{\perp}(d,N))$ agree (up to equivalence) with their classical counterparts.
\end{proposition}
\theproof
The only nontrivial case is ${\rm Th}(\mathbfcal{U}(d))$. But this follows immediately from the fact that the geometric realization preserves homotopy quotients along with naturality of the unit of the Quillen equivalence ${\rm sing}\dashv \vert \cdot \vert$.
\endofproof
 
\subsection{A smooth motivic spectrum model for Thom stacks}
\label{Sec mot}

In this section, we introduce the notion of an $D^1/\partial D^1$-spectrum and show that our 
geometrically refined cobordism stacks fit nicely into this setting. We also discuss some of the 
basic properties of these objects, leaving a more comprehensive discussion for a separate treatment.

\begin{definition}
[Motivic $D^1/\partial D^1$-spectra]
A {\rm $D^1/\partial D^1$-spectrum} $X$ is a sequence of pointed smooth stacks $X(n)$, 
equipped with structure maps 
$$
\xymatrix{
\sigma:D^1/\partial D^1\wedge X(n)
\ar[r] &
 X(n+1)\;.
}
$$
A morphism of $D^1/\partial D^1$-spectra is a sequence of level-wise maps $X(n)\to Y(n)$ commuting 
with the structure maps.
\end{definition}

The smash product and wedge product of pointed stacks are defined analogously to those of pointed spaces. 
That is, for two pointed stacks $X$ and $Y$, we define the {\rm wedge product} $X\vee Y$ as the pushout
$$
\xymatrix@R=1.3em{
\ast\ar[rr]\ar[d] && X\ar[d]
\\
Y \ar[rr] && X\vee Y
}
$$
and the {\rm smash product} as the quotient stack
$X\wedge Y:=X\times Y/(X\vee Y)$.

\medskip
The following examples will be of main interest.
\begin{example}[Madsen-Tillman spectrum]\label{mt spec}
Notice that the ordinary Madsen-Tillman spectrum ${\rm MT}(d)$ already defines a $D^1/\partial D^1$-spectrum via the usual maps
\(\label{usual MT}
\xymatrix{
D^1/\partial D^1\wedge {\rm Th}(\mathcal{U}^{\perp}(d,N))
\; \ar[r] & \;
 {\rm Th}(\mathcal{U}^{\perp}(d,N))
 }\;,
\)
with the subtle difference that we are taking quotients in a \emph{different} ambient category, namely, that of smooth stacks instead of topological spaces.
\end{example}
Notice that since we have an equivalence $q:\mathbfcal{U}^{\perp}(d,N)
\overset{\simeq}{\longrightarrow} \mathcal{U}^{\perp}(d,N)$, and this equivalence commutes with the structure maps \eqref{usual MT}, we can define the equivalent $D^1/\partial D^1$-spectrum as follows.

\begin{example}[Stacky Madsen-Tillman spectrum]\label{mt stack}
Let ${\rm Th}(\mathbfcal{U}^{\perp}(d,N))$ be the Thom stack of the universal orthogonal complement bundle 
over ${\bf Gr}(d,N)$. By Proposition \ref{pullback trivializes}, the pullback of the bundle 
$\mathbfcal{U}^{\perp}(d,N+1)\to {\bf Gr}(d,N+1)$ by the inclusion 
${\bf Gr}(d,N)\into {\bf Gr}(d,N+1)$ decomposes as $\mathbfcal{U}^{\perp}(d,N)\oplus {\bf 1}$. 
By Proposition \ref{Pro pro}, the resulting pullback map 
$\mathbfcal{U}^{\perp}(d,N)\oplus {\bf 1}\to \mathbfcal{U}^{\perp}(d,N+1)$ induces a map of 
Thom stacks
$$
\xymatrix{
D^1/\partial D^1\wedge{\rm Th}(\mathbfcal{U}^{\perp}(d,N))
\;\ar[r] & \;
{\rm Th}\big(\mathbfcal{U}^{\perp}(d,N+1)\big)
}\;.
$$
With these maps, ${\rm Th}(\mathbfcal{U}^{\perp}(d,N))$ becomes a smooth 
$D^1/\partial D^1$-spectrum with $d$ fixed and $N$ indexing the 
levels of the spectrum. We denote this spectrum by $\MM \TT(d)$. 
\end{example}

\vspace{1mm}
\begin{example} [Smooth motivic Thom spectrum] \label{tspec stack}
Let $\mathbfcal{U}(d)$ be the universal bundle over the classifying stack
$\BB {\rm O}(d)$. The pullback of $\mathbfcal{U}(d+1)$ to $\mathbf{B}{\rm O}(d)$ decomposes 
as $\mathbfcal{U}(d)\oplus {\bf 1}$ and we have an induced map at the level of Thom stacks
$$
\xymatrix{
D^1/\partial D^1\wedge {\rm Th}(\mathbfcal{U}(d))
\; \ar[r] & \;
{\rm Th}\big(\mathbfcal{U}(d+1)\big)
}\;.
$$
Thus, ${\rm Th}(\mathbfcal{U}(d))$ forms a smooth spectrum as $d$ varies. We denote this spectrum by 
$\mathbf{M}{\rm O}$. It is the smooth motivic model for the Thom spectrum.
\end{example}

\begin{remark}\label{st mt vs mt}
[Higher simplicial data]
{\bf (i)} We emphasize that there is no higher nondegenerate simplicial data in Example \ref{mt stack} and that we have an levelwise equivalence of $D^1/\partial D^1$-spectra
$${\rm MT}(d)\simeq \MM\TT(d)\;.$$ 

\medskip
\noindent {\bf (ii)} In contrast, Example \ref{tspec stack} has nondegenerate simplicial data in degree 1. We could also consider the the usual cobordism spectrum ${\rm MO}$ as a $D^1/\partial D^1$-spectrum (similar to the Madsen-Tillman spectrum), but this spectrum would not 
be equivalent to $\MM\mathbf{O}$.
\end{remark}

\begin{definition}
[$D^1/\partial D^1$ suspension- and infinite loop stacks]
\label{susloop}
{\bf (i)} Given a pointed stack $X$, we can define its \emph{$D^1/\partial D^1$-suspension spectrum} by successively smashing 
with $D^1/\partial D^1$, with structure maps being  identity. We denote this $D^1/\partial D^1$-spectrum 
by $\Sigma^{\infty}_{{}_{D^1/\partial D^1}}X$. 

\noindent {\bf (ii)} Conversely, if we are given a $D^1/\partial D^1$-spectrum $X(n)$, 
we can define its \emph{$D^1/\partial D^1$-infinite loop stack} as the colimit
$$
\Omega^{\infty}_{{}_{D^1/\partial D_1}}X(n):=
\colim_{n}\mathbf{Map}_{+}\big(D^n/\partial D^n,X(n)\big)\;,
$$
where the right hand side is the pointed mapping stack, defined in the obvious way.
\end{definition} 

\begin{remark}
[Motivic sphere vs. smooth sphere]
\label{smooth spheres}
Note that the sheaf $D^n/\partial D^n$ is not isomorphic to the smooth sphere $S^n$, viewed as a smooth manifold 
of dimension $n$. Indeed, if this were true, then an isomorphism would induce an isomorphism
$$
\hom(D^n/\partial D^n,\RR)\cong \hom(S^n,\RR)\;.
$$
Since the $\hom$ functor sends colimits to limits in the first variable, the Yoneda Lemma implies that 
the left side is in bijection the set of smooth functions on the closed disc (viewed as a manifold with 
boundary) $f:D^n\to \RR$, which restrict to a constant function on the boundary. Again by the Yoneda 
Lemma, the right side is in bijection with all smooth functions $f:S^n\to \RR$. But there are more 
elements in the former set, since any radially symmetric smooth function on the closed disc, with 
non-vanishing derivative at the boundary, will not define a smooth function on the sphere.
\end{remark}

\begin{proposition}
[Geometric realization of  $D^1/\partial D^1$-spectra]
\label{geometric realization of spectrum}
Let $X(n)$ be a $D^1/\partial D^1$-spectrum. The geometric realization of $X(n)$ is an ordinary 
(sequential, pre-) spectrum. Moreover, the $D^1/\partial D^1$-suspensions and infinite loop stacks 
geometrically realize to the usual suspension spectrum and infinite loop space.
\end{proposition}
\theproof
We first observe that since geometric realization preserves homotopy quotients (and in fact, strict quotients of sheaves), naturality of the unit of the adjunction $\vert \cdot \vert\dashv {\rm sing}$ implies that $\vert D^1/\partial D^1\vert\simeq \vert D^1 \vert/\vert \partial D^1 \vert\simeq S^1$. As a consequence, the geometric realization of the maps
$
\sigma: (D^1/\partial D^1)\wedge X(n)\longrightarrow X(n+1)
$
induce maps
$$
\xymatrix{
\vert \sigma \vert: S^1 \wedge \vert \Pi(X(n+1))\vert \ar[r]&
 \vert \Pi(X(n))\vert
}\;,
$$
 where $S^1$ is the topological circle. This proves the first claim. For the  geometric realization 
 of the $D^1/\partial D^1$-infinite loop space, the claim follows immediately from  from 
 Proposition \ref{eor pro}.
\endofproof

Proposition \ref{geometric realization of spectrum} implies that for the $D^1/\partial D^1$-spectra $\MM \TT(d)$, the $D^1/\partial D^1$-infinite loop stack (zero-truncated in this case) geometrically realizes to recover the topological 
infinite loop space of ${\rm MT}(d)$. The first of our two main theorems in the next section
 is devoted to refining the equivalence
\(
\label{GMTW}
\alpha: {\rm B}\mathscr{C}{\rm ob}\overset{\simeq}{\longrightarrow}
 \Omega^{\infty-1}{\rm MT}(d)\;,
\)
defined in \cite{GMTW}, to the smooth setting.

\begin{remark}
[Identifying the right stacky infinite loop space of the MT-spectrum]
 A natural choice for smooth analogue of the right hand side of \eqref{GMTW}
 would be the smooth stack $\Omega_{D^1/\partial D^1}^{\infty-1}\MM \TT(d)$ and, indeed, its 
geometric realization recovers the usual space on the right side of the above equivalence. However, it 
turns out that this is not quite right, and for crucial reasons:

\item {\bf (i)} The Pontrjagin-Thom construction translates cobordisms between smooth manifolds 
into \emph{smooth} paths between maps. The stack $\Omega_{D^1/\partial D^1}^{\infty-1}\MM \TT(d)$ does not contain the data of smooth paths between maps 
$$D^{d-1+N}/\partial D^{d-1+N}\to {\rm Th}(\mathbfcal{U}^{\perp}(d,N))\;.$$
This is extra data which needs to be accounted for. 

\item {\bf (ii)} Moreover, if we are to think of bordisms between smooth manifolds 
as time evolution of a smooth manifold (as motivated by physics and reflected by the definition in 
\cite{GMTW}), then the ordering on the interval $[0,1]$ (which indexes time) plays a crucial role. 
\end{remark}

All this information needs to be taken into account and motivates the following construction.

\paragraph{Time-ordering of maps.}
Fix a real number $0 < \epsilon <\!\!<1$. 
Let $X$ be a smooth stack and let $(t_0-\epsilon,t_1+\epsilon)\subset \RR$ be the open 
interval containing the closed interval $[t_0,t_1]\subset \RR$. One entity of main interest
will be $\mathbf{Map}^{\rm col}((\epsilon-t_0,t_1+\epsilon),X )$, the full substack of 
$\mathbf{Map}((\epsilon-t_0,t_1+\epsilon),X )$ on the maps 
$f:(t_0-\epsilon, t_1+\epsilon)\to X$ which make the following diagram commutative
$$
\xymatrix@R=1.5em{
(t_0-\epsilon,t_0] \ar[d] ~\ar@{^{(}->}[rr] && 
(t_0-\epsilon,t_1+\epsilon)\ar[d]^-{f} & & 
\ar@{_{(}->}[ll]~[t_1,t_1+\epsilon) \ar[d]
\\
\ast \ar[rr]^{{\rm ev_0}(f)} && 
X 
&& \ast \ar[ll]_{{\rm ev_1}(f)} 
}\;.
$$
This is just a fancy way of saying that restricting $f$ to the $\epsilon$-collars gives a constant map. 
For each $\epsilon$, we have a well-defined composition law
$$
m_{\epsilon}:\mathbf{Map}^{\rm col}
\big((\epsilon-t_0,t_1+\epsilon), X \big)\times_{X}
 \mathbf{Map}^{\rm col}\big((\epsilon-t_1,t_2+\epsilon), X\big )
\longrightarrow \mathbf{Map}^{\rm col}\big((\epsilon-t_0,t_2+\epsilon), X \big)\;,
$$
which on vertices concatenates smooth paths using the constancy on collars to glue. We get a directed
 system with respect to the ordering $(\epsilon,<)$ and the maps $m_{\epsilon}$ are compatible with the ordering. Thus, letting $\epsilon\to 0$ gives a composition in the limit. By abuse of notation, we will 
 denote the limit 
 $\mathbf{Map}^{\rm col}([t_0,t_1],X ):=\lim_{\epsilon\to 0}\mathbf{Map}^{\rm col}((t_0-\epsilon,t_1+\epsilon),X )$ and the limiting composition map by $m:=\lim_{\epsilon}m_{\epsilon}$. With this 
 composition map, we have an internal category given by the diagram
\(\label{concordance}
\xymatrix{
\RR \times X\sqcup \RR^2_{+}\times \mathbf{Map}^{\rm col}([0,1],X )
\; \ar@<.1cm>[r] \ar@<-.1cm>[r]& \; \RR \times X
}\;,
\)
where $\RR^2_{+}=C^{\infty}(-;\RR^2_+)$ is the sheaf of pairs of smooth functions $t_0$ and $t_1$ such that $t_1> t_0$ for all points in the domain. The source map evaluates at $\{0\}$ and picks out the map $t_0$ in the pair $(t_0,t_1)$. The target map evaluates at $\{1\}$ and picks out the function $t_1$. In particular, for the point Cartesian space $\ast$, we identify the factor $\{t_1,t_0\}\times {\bf Map}^{\rm col}([0,1],X)\cong {\bf Map}^{\rm col}([t_0,t_1],X)$, where the identification precomposes with the affine diffeomorphism $\varphi:(-\epsilon,1+\epsilon)\to (t_0-\epsilon,t_1+\epsilon)$. The composition map in \eqref{concordance} is a parametrized version of this composition. The source and target maps are induced by the usual evaluation maps.
\begin{definition}
[Smooth concordance category]
\label{def scc}
Let $X$ be a smooth stack. We define the 
\emph{smooth concordance category with ordering} ${\rm Conc}^{>}(X)$ 
to be the diagram \eqref{concordance} with composition map $m$ defined by the limiting 
operation on collarings above.\footnote{The superscript $>$ is used in our notation since the 
definition of this smooth category depends heavily on the standard ordering of $[0,1]$.}
 We denote the realization of the nerve of this groupoid as $\BB{\rm Conc}^{>}(X)$.
\end{definition}

Note that this new construct still geometrically realizes to the right object.

\begin{proposition}
[Topological correspondence] 
\label{Prop SHG}
Let $X\in \sh_{\infty}(\cartsp)$ be a smooth stack. We have a homotopy equivalence of geometric 
realizations
$$
\vert \BB {\rm Conc}^{>}(X)\vert \simeq \vert X \vert\;.
$$
\end{proposition}
\theproof
By the result of \cite{Pa} in Remark \ref{singular realization}, 
we have that the geometric realization of the stack $\mathbf{Map}^{\rm col}([0,1],X )$ 
can be computed as the simplicial set with $n$-vertices given by the colimit,
as $\epsilon\to 0$, of collared maps
$
(-\epsilon,1+\epsilon)\times \Delta^n\to X_n
$. 
Using the Cartesian closed structure on presheaves, we see that this simplicial set can be identified with 
$$
\xymatrix{
i:\map^{\rm col}\big((-\epsilon,1+\epsilon),{\bf sing}(X)\big) 
~\ar@{^{(}->}[r] &
~\map\big((-\epsilon,1+\epsilon),{\bf sing}(X)\big)
}\;,
$$
where the inclusion is that of a full sub-$\infty$-groupoid on the collared maps. Any such inclusion will induce 
an injection on connected components. By homotopy invariance of ${\bf sing}$, the projection $(-\epsilon,1+\epsilon)\to \ast$ induces an equivalence $j:{\rm sing}(X)\to \map((-\epsilon,1+\epsilon),{\bf sing}(X))$. Moreover, it is clear that the restriction of the map $i$ to the components of the constant maps $(-\epsilon,1+\epsilon)\to\ast \to  {\bf sing}(X)$ factors through $j$ and is surjective. Thus, $i$ induces an isomorphism on connected components. By definition, $i$ induces an isomorphism on higher homotopy groups and defines an equivalence of Kan complexes. Since $\epsilon>0$ was arbitrary, the homotopy colimit as $\epsilon\to 0$ is the homotopy colimit over a constant diagram and is therefore equivalent to ${\rm sing}(X)$. Since strict filtered colimits in $\sh_{\infty}(\cartsp)$ model their homotopy colimits, we conclude that 
$$
\map^{\rm col}\big([0,1],{\bf sing}(X)\big)
\simeq {\rm sing}(X) \simeq \Pi(X)\;.
$$
From this calculation we see that geometrically realizing the nerve of the internal category \eqref{concordance} gives the simplicial diagram in CGWH spaces
$$
\Big\{\xymatrix{
\hdots \ar@<.1cm>[r]\ar[r]  \ar@<-.1cm>[r] & \RR\times \vert X\vert\sqcup \RR^2_+\times \vert X \vert
\; \ar@<.05cm>[r] \ar@<-.05cm>[r]& \; \RR\times \vert X\vert \;
}\Big\}\cong {\cal N}(\mathscr{C})\times \vert X\vert\;,
$$
where ${\cal N}(\mathscr{C})$ is the nerve of the category with a single object and two morphisms $1$ and $a$ satisfying $a^2=a$. \footnote{This identification is obtained by projecting out both $\RR$ and $\RR^2_+$, which is levelwise a weak equivalence commuting with the face maps.} But the Kan fibrant replacement of ${\cal N}(\mathscr{C})$ is clearly contractible. Therefore, the geometric realization of the simplicial space on the right is equivalent to 
$\vert X \vert$. Finally, the claim follows from Proposition \ref{Prop real}.
\endofproof

From Proposition \ref{geometric realization of spectrum}, we immediately get the following 
crucial characterization.

\begin{corollary}
[Correspondence of infinite loop spaces]
\label{realization of concordance}
The geometric realization of the smooth stack\;
 $\BB{\rm Conc}^{>}\big(\Omega^{\infty-1}_{{}_{D^1/\partial D^1}}{\MM \TT}(d)\big)$
is equivalent to the infinite loop space $\Omega^{\infty-1}{\rm MT}(d)$.
\end{corollary}

\medskip
In the abstract Pontrjagin-Thom construction (Section \ref{Sec PT}), we will need to make use of the stacky homotopy type of the smooth stack
$\BB{\rm Conc}^{>}\big(\Omega^{\infty-1}_{{}_{D^1/\partial D^1}}{\MM\TT}(d)\big)$. 
In high degrees, this turns out to be trivial. 

\begin{proposition}
[Sheaves of homotopy groups]
\label{Prop SHG}
For each $t\in \RR$, indexing a map $t:\ast\to \RR$, we have a canonical inclusion
$$
\xymatrix{
t:X
\; \ar@{^{(}->}[r] & \;
\BB{\rm Conc}^{>}\big(X\big)
}\;,
$$
corresponding to the inclusion $t:X\into \RR\times X$. For each $t$, this inclusion induces an isomorphism on $\widetilde{\pi}_k$, based in the $t$-component, for all $k\geq 1$.
\end{proposition}
\theproof
Consider the prestack which models $\BB{\rm Conc}^>(X)$. From the definition of the realization, we see that this prestack can be modeled by the diagonal of the bisimplicial diagram given by taking the nerve of \eqref{concordance}. The first few levels look like
\(
\hspace{-.1cm}
\label{pre concordance}
\BB{\rm Conc}^>(X)=\Big\{\xymatrix@=1.5em{
\hdots \ar@<.1cm>[r]\ar[r]  \ar@<-.1cm>[r] & \RR \times X_1\sqcup  \RR^{2}_+ \times {\bf hom}^{\rm col}([0,1],X_1)
\; \ar@<.05cm>[r] \ar@<-.05cm>[r]& \; \RR\times X_0
}\Big\}
\)
where ${\bf hom}$ denotes the presheaf of maps. Consider the prestack given by
\(
\hspace{-.1cm}
\label{pre concordance}
Y=\Big\{\xymatrix@=1.5em{
\hdots \ar@<.1cm>[r]\ar[r]  \ar@<-.1cm>[r] & \RR \sqcup  \RR^{2}_+ 
\; \ar@<.05cm>[r] \ar@<-.05cm>[r]& \; \RR
}\Big\}\;.
\)
Then the projection map $p:\BB{\rm Conc}^>(X)\to Y$ is an objectwise Kan fibration. Moreover for any fixed map $x:\ast\to U$, the induced map ${\rm ev}_x:Y(U)\to Y(\ast)\simeq {\cal N}(\RR^{\delta},\leq)$ is the Kan fibration which simply evaluates the smooth functions at a fixed point. Thus for each $U\in \cartsp$, we have a Kan fibration
$$
\xymatrix{
 {\rm ev_x} \circ p :\BB{\rm Conc}^>(X)(U) \; \ar[r] & \;  {\cal N}(\RR^{\delta},\leq) 
 }\;,
$$
and the fiber at $r\in \RR$ can be identified with the subspace $X(U)\times S(U)\subset X(U)\times C^{\infty}(U,\RR)$ on those smooth functions $f:U\to \RR$ which satisfy $f(x)=r$. Since $C^{\infty}(U;\RR)$ is zero-truncated, the inclusion at the constant map $r:U\to \RR$, which maps every point in the domain to $r$ induces an isomorphism
$$\pi_k(X(U))\overset{\cong}{\longrightarrow} \pi_k(X(U)\times S(U))\;.$$ 
for $k\geq 1$. Since ${\cal N}(\RR^{\delta},\leq )$ is contractible, we see that ${\rm ev_x} \circ p$ must induce an isomorphism on $\pi_k$ for $k\geq 1$.
\endofproof

Since the stack $\Omega^{\infty-1}_{{}_{D^1/\partial D^1}}{\MM \TT}(d)$ is zero-truncated, we immediately  deduce the following.
\begin{proposition}
The smooth stack $
\BB{\rm Conc}^{>}\big(\Omega^{\infty-1}_{{}_{D^1/\partial D^1}}{\MM \TT}(d)\big)
$
is zero-truncated.
\end{proposition}

\section{A stacky perspective on the cobordism category}
\label{Sec stacky}

\subsection{The smooth cobordism category of Galatius-Madsen-Tillman-Weiss}
\label{Sec GMTW}

In this section we describe a variant of the topological cobordism category, where we regard 
both the space of objects and morphisms as smooth objects. Our definition will be closely related 
to the sheaf of categories describing the parametrized cobordisms introduced in 
\cite{GMTW}. In fact, if we work instead over the site of smooth manifolds, we will recover 
exactly the sheaf of categories described there. 

\begin{definition}
[Cobordism category \cite{GMTW}]
\label{framed cobordisms}
The {\rm $d$-dimensional cobordism category}  $\mathscr{C}{\rm ob}_d$ 
 is the category with objects given by pairs $(M,t)$, 
with $M\subset \RR^{\infty-1}=\colim_N\RR^{N-1}$ a smooth submanifold and $t\in \RR$. 
The morphisms are given by triples $(W,t_0,t_1):(M_0,t_0)\to(M_1,t_1)$, such that 

\vspace{.1cm}

\noindent {\bf (i)} {\rm (Ordering)} $t_0<t_1$.
\\
\noindent {\bf (ii)} {\rm (Embedding)}  $W\subset \RR^{\infty-1}\times [t_0,t_1]$.
\\
\noindent {\bf (iii)} {\rm (Collared neighborhood)} There is $\epsilon>0$ such that $W$ 
restricted to $[t_0,t_0+\epsilon)$ is 
$M_0\times [t_0, t_0+\epsilon)$ and the restriction of $W$ to $(t_1-\epsilon,t_1]$ is 
$M_1\times (t_1-\epsilon,t_1]$.
\\
\noindent {\bf (iv)} {\rm (Boundary condition)} 
$W\cap \{t_0,t_1\}\times \RR^{\infty-1}=\partial W$.
\end{definition}
This definition contains quite a bit of data. However, the complexity of the definition 
is necessary to ensure that one has well-defined compositions of morphisms.
\begin{figure}[h]\label{bordism composition}
\centering

\begin{tikzpicture}[scale=.7]

\draw[dashed] (-3,4) -- (-3,-2) -- (11,-2) -- (11,4) -- (-3,4) node[at end, below right] {$\RR^n$};

\draw (-1,0) .. controls (-1,1) and (-2,1) .. (-2,2) node[near start, left] {$M$};
\draw[dashed] (-.2,.35) .. controls (-.2,1.35) and (-1.2,1.35) .. (-1.2,2.35);
\draw (-2,2)..controls (0,3) and (1,2)..(4,2) node[midway, above] {$W$}; 
\draw (-1,0)..controls (1,1) and (2,0)..(4,0); 
\draw (4,0)..controls (4,.5) and (3.5,1.5)..(4,2) node[midway, right] {$N$};
\begin{scope}[yscale=1,xscale=-1,xshift=-6cm]
\draw[dashed] (3,.2)..controls (3,.5) and (3.5,1.5)..(3,2);
\end{scope}

\begin{scope}[yscale=1,xscale=-1,xshift=-8cm]
\draw (-1,0) .. controls (-1,1) and (-2,1) .. (-2,2)  node[near start, right] {$M^{\prime}$};
\draw[dashed] (-.2,.35) .. controls (-.2,1.35) and (-1.2,1.35) .. (-1.2,2.35);
\draw[arrows=<->] (-1,1.62)--(-1.8,1.35) node[midway,above] {$\epsilon$};
\draw (-2,2)..controls (0,3) and (1,2)..(4,2) node[midway, above] {$W^{\prime}$}; 
\draw (-1,0)..controls (1,1) and (2,0)..(4,0); 
\draw[dashed] (3,.2)..controls (3,.5) and (3.5,1.5)..(3,2);
\end{scope}

\end{tikzpicture}
\caption{A composition of cobordisms with collared neighborhoods.}
\end{figure}
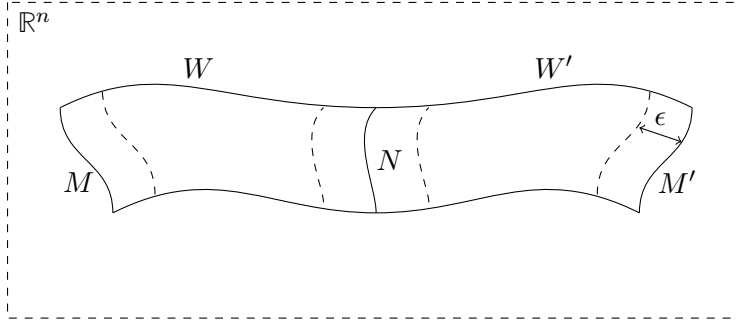

This category admits a smooth structure in a natural way, i.e. it can be made into an internal category 
in smooth sheaves ${\Sh}(\cartsp)$. To see this, let us recall that if ${\rm Emb}(M,\RR^{d+\infty-1})$ 
denote the set of embeddings of closed $(d-1)$-dimensional manifolds and ${\rm Diff}(M)$ denote the 
diffeomorphism group of $M$, then both of these spaces admit infinite-dimensional, smooth, Fr\'echet 
structure (see for instance \cite{KM}). Moreover, the quotient space 
${\rm Emb}(M,\RR^{d+\infty-1})/{\rm Diff}(M)$, where diffeomorphisms act freely by 
precomposition, also admits the structure of a smooth infinite-dimensional manifold. In \cite{BF}\cite{KM}, 
it is shown that ${\rm Emb}(M,\RR^{d-1+\infty})/{\rm Diff}(M)$ is the base space of a smooth fiber bundle 
$$
\xymatrix{
\pi:{\rm Emb}(M,\RR^{d-1+\infty})\times_{{\rm Diff}(M)}M
\; \ar[r] & \;
{\rm Emb}(M,\RR^{d-1+\infty})/{\rm Diff}(M)
}\;,
$$
and that this bundle is a universal bundle which describes smooth manifolds parametrized over 
a base manifold $B$. More precisely, given a map 
$$
\xymatrix{
f:B\ar[r] & {\rm Emb}(M,\RR^{d+\infty-1})/{\rm Diff}(M)\;,
}
$$
the pullback of $\pi$ by $f$ gives a submanifold $E\subset B\times \RR^{d-1+\infty}$ and the projection 
to $B$ defines a smooth fiber bundle with fiber $M$.

\medskip
Given that we have a smooth structure available, we can regard each of these objects as a 
smooth sheaf on the site of Cartesian spaces, by embedding them via the sheaf of smooth plots. 
For bordisms between manifolds, we can consider the analogous smooth object constructed as follows. 
Let $W$ be a compact $d$-manifold with collared boundary of width $\epsilon>0$ (as depicted in Figure 1 above)  and let ${\rm Emb}_{\epsilon}(W,[0,1]\times \RR^{d+\infty-1})$ be the smooth space of embeddings of $W$. Define the infinite-dimensional smooth manifold 
$$
{\rm Emb}(W,~[0,1]\times \RR^{d+\infty-1}):=
\colim_{\epsilon>0}{\rm Emb}_{\epsilon}(W,~[0,1]\times \RR^{d+\infty-1})\;,
$$
where the colimit is taken over the obvious maps as $\epsilon\to 0$.  Let ${\rm Diff}_{\epsilon}(W,\partial_{\rm in},\partial_{\rm out})$ denote the subgroup of the diffeomorphism group of $W$ which restricts on the collars to diffeomorphisms of the form $\phi\times {\rm id}$, with $\phi$ a diffeomorphism of the boundary. Similar to the embeddings, we set
$$
{\rm Diff}(W):=\colim_{\epsilon>0}{\rm Diff}_{\epsilon}(W,\partial_{\rm in},\partial_{\rm out})\;.
$$
Again, we can divide out the ${\rm Diff}(W)$-action on the embeddings and we get a smooth fiber bundle which classifies parametrized bordisms. This leads us to consider the following. 
\begin{definition}
[Smooth cobordism category]
\label{scc}
The \emph{smooth cobordism category $\mathscr{C}{\rm ob}_d$} is the category with 
sheaf of objects and stack of morphisms given, respectively, by
\begin{align*}
{\rm Ob}(\mathscr{C}{\rm ob}_d) &:=
\underset{M}{\coprod}\RR\times {\rm Emb}(M,\RR^{d+\infty-1})/{\rm Diff}(M)\;,
\\
{\rm Mor}(\mathscr{C}{\rm ob}_d) & := {\rm Ob}(\mathscr{C}{\rm ob}_d) \sqcup
\underset{W}{\coprod}\RR^2_+\times {\rm Emb}(W,\;[0,1]\times \RR^{d+\infty-1})/{\rm Diff}(W)\;.
\end{align*}
\end{definition}
Under the geometric realization $\vert \cdot \vert \circ \Pi: {\Sh}(\cartsp) \overset{\Pi}{\longrightarrow} \sset 
\overset{\vert \;\cdot \; \vert}{\longrightarrow} \mathscr{T}{\rm op}$, 
this smooth category recovers the topological cobordism category. This follows from 
two facts: First, that the geometric realization is a left adjoint at the level of 1-categories 
and hence preserves quotients; second, that the geometric realization of a smooth 
manifold is equivalent to the underlying topological space of the manifold.

\begin{remark} 
[Strict vs. homotopy quotient]
Note that we have used a strict quotient in our definition of the smooth cobordism category. 
In fact, since ${\rm Diff}(M)$ acts freely and transitively on 
${\rm Emb}(M,\RR^{d+\infty-1})$, the homotopy quotient is modeled by the strict quotient. Thus, we could have (more generally) defined the smooth category via the homotopy quotient. However, the particular model we chose will be more useful in calculation. 
\end{remark}

Note that, being an internal category in the category of smooth sheaves, we can regard 
$\mathscr{C}{\rm ob}_d$ as a sheaf of categories on the site of smooth manifolds, 
equipped with the usual topology of open covers, by sending each smooth manifold 
$M\in \mathscr{M}{\rm an}\into \Sh(\cartsp)$ to the category with objects 
$\hom(M,{\rm Ob}(\mathscr{C}{\rm ob}_d))$ and morphisms 
$\hom(M,{\rm Mor}(\mathscr{C}{\rm ob}_d))$. 
This sheaf of categories was considered 
in \cite{GMTW}.

\medskip
How can we think of the objects of the new category?
Notice that, for each $U\in \cartsp$ an element in the set of objects ${\rm Ob}(\mathscr{C}{\rm ob}_d)(U)$ 
is simply a trivializable bundle $N\cong M\times U \to U$ with fiber a $(d-1)$-dimensional $M$. Moreover, each 
$N\subset U\times \RR^{d-1+\infty}$ and the projection onto $U$ gives the bundle map. The bundle $N$ comes equipped with a smooth function $t:U\to \RR$, indexing the position of the fibers in time. The triviality is a consequence of the fact that we are working over the small site of convex subsets 
of Euclidean spaces, where all smooth fiber bundles trivialize. Similarly, one gets an identification 
of the set of morphisms (after evaluating at $U$) as trivial bundles $Z\cong W\times U \to U$, 
with boundary $\partial Z\simeq N\sqcup N^{\prime}$ and a pair of smooth functions $t,t^{\prime}:U\to \RR^2$, for which $t^{\prime}(u)>t(u)$ for all $u\in U$. In \cite{GMTW} a slightly different sheaf of categories, 
isomorphic to $\mathscr{C}{\rm ob}_d$, is introduced and this category is defined on the site 
of smooth manifolds. Since the resulting bundles may be nontrivial, their sheaf of 
categories involves a slightly more complicated definition. 

\begin{definition}[Classifying stack of the cobordism category]
\label{stacky cobordism category}
We define the \emph{smooth cobordism stack $\BB\mathscr{C}{\rm ob}_d$} 
as the stacky realization of the smooth cobordism category $\mathscr{C}{\rm ob}_d$
(from Definition \ref{scc}).
\end{definition}
Under geometric realization, this stack indeed recovers the usual cobordism category.
\begin{proposition}
[Classical correspondence]
\label{prop TBC}
The geometric realization of $\BB\mathscr{C}{\rm ob}_d$ is weak equivalent to 
the classifying space of the topological bordism category $B\mathscr{C}{\rm ob}_d$.
\end{proposition}
\theproof
This follows immediately from Proposition \ref{Prop real}.
\endofproof

As in the case of the concordance category for the Madsen-Tillman spectrum, we also have the following.
\begin{proposition}
[Stacky triviality] 
\label{BBC0}
The smooth stack $\BB\mathscr{C}{\rm ob}_d$ is zero-truncated.
\end{proposition}
\theproof
The same proof as that of Proposition \ref{Prop SHG} works here as well.
\endofproof

\subsection{An abstract Pontrjagin-Thom construction}
\label{Sec PT}

In this section, we prove the first of our main theorems  via a smooth refinement of the Pontrjagin-Thom construction. This construction will take place entirely in the category of smooth stacks. It will also refine the usual Pontrjagin-Thom collapse map in the sense that, if we geometrically realize our stacks, we will recover the usual collapse map.

\medskip
Before presenting the construction of the Pontrjagin-Thom map and 
the first of our main theorems, we will need a few lemmas which will 
be used in the construction and the proof. The first gives us a nice 
characterization of the hypercovers of the 
sphere $D^{d+N}/\partial D^{d+N}$, the second will provide a smooth family of diffeomorphisms, while the third connects these to regular values at zero. 

\begin{lemma}
[Covers for $D^n/\partial D^n$]
\label{Lem cover}
Consider the smooth unit sphere $S^{d+N}\subset \RR^{d+N+1}$, equipped with a choice of basepoint $\infty:\ast\to  S^{d+N}$. Let $\{U_{\beta}\}$ be a good open cover of $S^{d+N}-\{\infty\}$ and let $V$ is a sufficiently small geodesically convex open neighborhood of $\{\infty\}$, so that
$Y:=\coprod_{\beta}U_{\beta}\sqcup V$
is a good open cover of $S^{d+N}$. Let 
$u:D^{d+N}/\partial D^{d+N}\to S^{d+N}$ 
be the map induced by the universal property of the quotient. Then 

\item {\bf (i)} the pullback $Y^{\prime}:=u^*Y$ is defines an effective epimorphism $Y\to D^{d+N}/\partial D^{d+N}$. 

\item {\bf (ii)} Moreover, the pullback is of the form 
$
Y^{\prime}=\coprod_{\beta}W_{\beta}\sqcup Z
$,
with $\{W_{\beta}\}$ a cover of the interior $(D^{d+N})^o$ (in the traditional sense) and 
$Z$ the pullback of $V$. 
\end{lemma}
\theproof
The first claim is immediate, since effective epimorphisms are stable under pullback. The second claim follows from the fact that coproducts commute with pullbacks in any topos (i.e. coproducts are universal).
\endofproof

\vspace{1mm}
\begin{lemma} [A smooth family of diffeomorphisms]
\label{lem sfd}
Fix $1 > \epsilon>0$. Then, for each real number $0<\delta<\epsilon$, and for any point
$y\in D^N_{\delta}\subset D^N$, there is a smooth family of diffeomorphisms 
$\varphi^y_t:D^N\to D^N$, $t\in[0,1] $, such that 
\begin{enumerate}[label={\bf (\roman*)}]

\vspace{-2mm}
\item $\varphi^y_0={\rm id}$, $\varphi^y_1(y)=0\in D^N$, and 

\vspace{-2mm}
\item for each pair of equidistant points $y,y^{\prime}\in D^N_{\delta}$, 
there is fixed orthogonal transformation $Q^{y,y^{\prime}}:D^N\to D^N$ satisfying
$
\varphi^{y^{\prime}}_t\circ Q^{y,y^{\prime}}=Q^{y,y^{\prime}}\circ \varphi_t^{y}
$,
for all $t$.
\end{enumerate}
\end{lemma}
\theproof
Let $0<\delta<\epsilon$. Fix a radially symmetric smooth bump function $\rho_{\epsilon}$ on $D^N$ which has support in the interior of $D^N_{\epsilon}$ and which is 1 on $D^N_{\delta}$. Let $y\in D^N_{\delta}$ and consider the smooth family of diffeomorphisms
$$
\varphi^y_t(x):=x-\rho_{\epsilon}(x)ty\;.
$$
Then $\varphi^y_0(x)=x$, $\varphi^y_1(y)=y-\rho_{\epsilon}(y)y=0$. Moreover if $y$ and $y^{\prime}$ are equidistant, we can choose an orthogonal transformation $Q^{y,y^{\prime}}$ taking $y$ to $y^{\prime}$. Then 
$$
Q^{y,y^{\prime}}\circ \varphi^{y}_t=Q^{y,y^{\prime}}(x)+\rho_{\epsilon} ty^{\prime}=\varphi^{y^{\prime}}_t(Q^{y,y^{\prime}}(x))\;.
$$

\vspace{-10mm}
\endofproof

\begin{lemma}
[Regular values on the disc]
 \label{rvd}
Fix $1 > \delta>0$. Let $\{f_k\}_{k=1}^{\infty}:D^N\to D^N$ be any countable collection of smooth functions. Then there is disc of radius $0\leq \gamma\leq \delta$ whose boundary contains regular values 
of all the $f_k$'s. \footnote{We use the convention that the boundary of the trivial disc $\{0\}$ is again 
the point $\{0\}$.}
\end{lemma}
\theproof
Suppose that is not the case. Then for all discs of radius $D^N_{\gamma}\subset D^N_{\delta}$, there is a function $f_k\in \{f_k\}_{k=1}^{\infty}$ such that the boundary contains no regular values of $f_k$. Let $S_k$ be the set of all discs (including the degenerate case $\{0\}$) with no regular values of $f_k$ on the boundary. Let $E_k\subset D^N_{\delta}$ be the union of all the boundaries of the discs in $S_k$. Since $D^N_{\delta}$ is the union of all $N$-spheres centered at zero of arbitrary radius $0\leq \gamma <\delta$, by hypothesis, we can write
$$
D^N_{\delta}=\bigcup_{k=1}^{\infty}E_k\;.
$$  
Let $C_k$ denote the set of all critical values of $f_k$. By Sard's theorem, each $C_k$ has measure zero in $D^N$. If each $E_k\subset C_k$ was measurable, then the right hand side would be a countable union of measure zero sets, and hence would have measure zero. Clearly the left hand side has positive measure and we would arrive at the desired contradiction.

It remains to prove that each $E_k$ is measurable and to this end, it suffices to prove that each $E_k$ is 
closed. Indeed, let $\{x_m\}$ be a sequence in $E_k$ and suppose $x_m\to x_{\infty}\in D^N_{\delta}$. Since the set of critical values $C_k$ of $f_k$ is closed, and every point on $E_k$ is a critical point, we have that $x_{\infty}\in C_k$. We claim that actually $x_{\infty}\in E_k$. Fix any point on the sphere of radius 
$\vert x_{\infty}\vert\geq 0$ and an orthogonal transformation $Q$ sending $x_{\infty}$ to $y$. 
Since $Q$ is an orthogonal transformation and $x_m\in E_k$ for all $m$, we must have $Q(x_m)\in E_k$ 
for all $m$. Since $Q$ is continuous on the disc, we have
$$
\lim_{m\to \infty}Q(x_{m})=Q(x_{\infty})=y\;.
$$
Since $C_k$ is closed we, therefore,
 have $y\in C_k$. Since $y$ was an arbitrary point on the sphere 
of radius $\vert x_{\infty}\vert$, we must have $x_{\infty}\in E_k$.  
\endofproof

\medskip
\noindent {\bf The collapse map.}
This map will be a morphism of smooth sheaves
\(\label{abstract collapse map}
\xymatrix{
{\bf PT}:\widetilde{\pi}_0\big(\BB\mathscr{C}{\rm ob}_{d}\big)
\ar[r] & \;
\widetilde{\pi}_0\big(\BB{\rm Conc}^{>}\big(\Omega^{\infty-1}_{{}_{D^1/\partial D^1}}{\MM\TT}(d)\big)\big)
}\;,
\)
where the source and target categories are defined in Section \ref{Sec GMTW} (see Definition \ref{scc}) and Section \ref{Sec mot} (see Definition \ref{def scc}, Corollary \ref{realization of concordance} and Proposition \ref{Prop SHG}),
respectively.  Recall moreover that since $\Omega^{\infty-1}_{D^1/\partial D^1}{\MM\TT}(d)$ is zero-truncated as a smooth stack and since $\MM\TT(d)$ is equivalence to the $D^1/\partial D^1$-spectra of strict sheaves ${\rm MT}(d)$ (see Remark \ref{st mt vs mt}), up to equivalence the smooth category ${\rm Conc}^{>}\big(\Omega^{\infty-1}_{{}_{D^1/\partial D^1}}{\MM\TT}(d)\big)$ has sheaf of objects and morphisms given, respectively, by
\begin{align*}
{\rm Ob} &:=\colim_{N}\RR\times {\bf hom}_{+}\big(D^{d+N-1}/\partial D^{d+N-1},\; {\rm Th}(\mathcal{U}^{\perp}(d,N))\big)\;,
\\
{\rm Mor} &:=\colim_N\RR^2_+\times {\bf hom}_{+}^{\rm col}\big([0,1]_+\wedge D^{d+N-1}/\partial D^{d+N-1}, \; {\rm Th}(\mathcal{U}^{\perp}(d,N))\big)\sqcup \ {\rm Ob}\;,
\end{align*}
where ${\bf hom}_+(-,-)$ denotes the sheaf of pointed maps. 

\medskip
Fix a Cartesian space $U\in \cartsp$. After evaluating the sheaf of morphisms at $U$, an element of the resulting set is a choice of submanifold $N\subset U\times \RR^{d+N-1}$, such that the map $\pi:N\to U$ is a trivializable bundle, with fiber $M$. We, therefore,
 consider the diagram 
\(\label{maps to universal}
\xymatrix{
T^{\perp}N\cong T^{\perp}M\times U 
\ar[rrr]^-{\gamma}\ar[d]_{{\rm id}\times \pi} &&&
\mathcal{U}^{\perp}(d,N)\ar[r]\ar[d] & \mathcal{U}^{\perp}(d)\ar[d]
\\
N\cong M\times U \ar[rrr]^-{T^{\pi}M} 
&&& {\rm Gr}(d,N)\ar[r] & {\rm Gr}(d,\infty)\;,
}
\)
where $T^{\pi}M$ is the fiberwise Gauss map, sending a pair $(x,u)\in M\times U$ to $T_xM\times \RR\subset \RR^{d+N}$. The map  $\gamma$ defined fiberwise as the canonical map which sends a pair $(v,u)\in T^{\perp}M\times U$, $v\in T^{\perp}_xM$, to $v\in \RR^N$. Taking a tubular neighborhood of $M$ in some disc with sufficiently large radius $D^{d+N-1}$ and collapsing the complement \footnote{Note that we are collapsing in smooth stacks and not in spaces.} gives a map
$$
\xymatrix{
(D^{d+N-1}/\partial D^{d+N-1})\wedge U_{+}
\ar[r] &  {\rm Th}(T^{\perp}M)\wedge U_+ \ar[r] &
{\rm Th}\big(\mathcal{U}^{\perp}(d,N)\big)\;,
}
$$
which is independent of the choices made up to (collared) concordance. Similarily, we can define the map at the level of bordisms, using the Gauss map sending $(y,u)\mapsto T_yW\subset \RR^{d+N}$, to get a (collared) concordance
$$
\xymatrix{
(D^{d+N-1}/\partial D^{d+N-1})\wedge [0,1]_+\wedge U_{+}
\ar[r] &  {\rm Th}(T^{\perp}W)\wedge U_+ \ar[r] &
{\rm Th}\big(\mathcal{U}^{\perp}(d,N)\big)\;,
}
$$
which depends on a choice of tubular neighborhood and the radius of the disc $D^{d+N-1}$. This shows that the map ${\bf PT}$ is well defined and natural in $U$.

\medskip
We now claim that the collapse map defined above induces a weak equivalence of smooth stacks. In the introduction we claimed that the proof of our main theorem has the advantage 
of avoiding some rather delicate differential topology theorems (in particular, Phillips' submersion 
theorem \cite{Ph}) used in the proof of \cite{GMTW}
(see \cite{Fr2} for an exposition). We will still need a bit of {\it local} transversality, 
\footnote{By local, we mean a condition that holds on local patches and is appropriately compatible on intersections.} and the price we pay for simple local conditions is complicated gluing conditions. 
 We now proceed with this setup.

\begin{remark} [Transversality] We could probably have avoided \emph{all} transversality issues 
had we chosen to work with derived manifolds and derived cobordisms between them. However, 
since we wanted to directly connect the cobordism category of \cite{GMTW} with our smooth 
motivic motivic model for the Madsen-Tillman spectrum, we chose to work with smooth manifolds
 instead. 
\end{remark}

Considering now the stacky model $\mathbfcal{U}^{\perp}(d,N)$ for $\mathcal{U}^{\perp}(d,N)$, recall that a morphism $f:M\to {\rm Th}(\mathbfcal{U}^{\perp}(d,N))$ is secretly defined as a map out of 
the {\v C}ech nerve $\check{C}(\{U_{\alpha}\})$ of some good open cover of $M$. Unravelling the resulting cocycle data, we see that such a map in particular gives an assignment morphisms 
$f_{\alpha}:U_{\alpha}\to D^{N}/\partial D^N$ on open sets, which on intersections obeys the compatibility relation $f_{\beta}=g_{\alpha\beta}f_{\alpha}$, for $g_{\alpha\beta}:U_{\alpha\beta}\to {\rm O}(N)$. 
We will say that the map $f$ is \emph{regular} if for all $\alpha$ such that 
$f_{\alpha}:U_{\alpha}\to D^N/\partial D^N$ factor through $D^N$, $0\in D^N$ is a regular 
value of the smooth function $f_{\alpha}$. 

\begin{proposition}
[A homotopy for the classifying space of concordance]
\label{homotopy edge}
Choose a point $f:\ast\to \BB{\rm Conc}^{>}(\Omega^{\infty}_{{}_{D^1/\partial D^1}}\MM \TT(d))(U)$. Such a point is equivalently a choice of map 
$$
\xymatrix{
f:D^{d-1+N}/\partial D^{d-1+N}\wedge U_+
\ar[r] &
 {\rm Th}(\mathbfcal{U}^{\perp}(d,N))
 }
$$
and a smooth function $a:U\to \RR$. Then there exists an edge
\(
\label{eq edges}
\xymatrix{
H:\Delta[1] 
\; \ar[r] & \;
\BB{\rm Conc}^{>}(\Omega^{\infty-1}_{{}_{D^1/\partial D^1}}\MM \TT(d))(U)
}\;,
\)
such that $d_0H=(f,a)$ and $d_1H=(g,a^{\prime})$, with $a^{\prime}>a$ and $g$ regular. \footnote{Here $d_0$ and $d_1$ denote the evaluations at the respective endpoints of the edge.} Moreover, $H$ can be chosen so that the preimages $N_{\beta}=g_{\beta}^{-1}(0)\cong M_{\beta}\times U$, with $g_{\beta}$ the restriction of $g$ to a specified cover of $D^{d-1+N}/\partial D^{d-1+N}$ and $M_{\beta}\subset D^{d-1+N}$ a submanifold.
\end{proposition}
\theproof
An edge \eqref{eq edges} is a choice of collared maps
$$
\xymatrix{
H:[0,1]_+\wedge D^{d-1+N}/\partial D^{d-1+N}\wedge \Delta[1]_+\wedge U_+ 
\; \ar[r] & \;
 {\rm Th}(\mathbfcal{U}^{\perp}(d,N)) 
 } \;,
$$
and a pair of smooth functions $a<a^{\prime}:U\to \RR$. The face maps are given by precomposing with the composition of coface maps $\delta^0,\delta^1:\{\ast\}\into [0,1]$ and $d^0,d^1:\Delta[0]\into \Delta[1]$. The face maps send the ordered pair $(a,a^{\prime})$ to $a$ and $a^{\prime}$, respectively. In particular, we can take the collared map to be of the form
$$
\xymatrix{
H:[0,1]_+\wedge D^{d-1+N}/\partial D^{d-1+N}\wedge U_+ 
\; \ar[r] & \;
 {\rm Th}(\mathbfcal{U}^{\perp}(d,N)) 
 } \;,
$$
by considering it as degenerate in the simplicial direction. 

We need to identify what such an edge like this looks like in terms of cocycle data. To this end, recall that
${\rm Th}(\mathbfcal{U}^{\perp}(d,N))$ can be identified with the smooth stack given 
level-wise by
$$
{\rm Th}(\mathbfcal{U}^{\perp}(d,N))=\left\{\xymatrix{
\hdots \ar@<.1cm>[r]\ar[r]\ar@<-.1cm>[r] & 
({\rm O}(N)\times {\rm V}(N,d))_+\wedge D^N/\partial D^N\ar@<.05cm>[r]\ar@<-.05cm>[r] &  
{\rm V}(N,d)_+\wedge D^N/\partial D^N
}\right\}\;.
$$
Let $\{W_{\beta}\}\cup \{Z\}$ be a generalized cover of $D^N/\partial D^N$ of the form presented in Lemma \ref{Lem cover}.  After 
passing to this cover and working out the cocycle data, we see that a map $H$ is 
uniquely determined by the following maps
%
%
 \begin{align*}
 & (a,a^{\prime}): U
 \longrightarrow 
 \RR^2_+
 \\
 H^{\beta} : \; &[0,1]\times U\times W_{\beta}
\longrightarrow 
D^{N}/\partial D^{N},
\qquad \qquad \quad
 g^{\alpha\beta}: [0,1]\times U\times  W_{\alpha\beta}
 \to {\rm O}(N)\times {\rm V}(N,d),
 \\
  H^{Z}: \;& [0,1]_+\wedge U_+\wedge Z
\longrightarrow 
D^{N}/\partial D^{N},
 \qquad \qquad \;
  g^{Z\beta}:[0,1]\times U\times  W_{\beta}\cap Z
 \to {\rm O}(N)\times {\rm V}(N,d),
\end{align*}
 which satisfy the the usual compatibility condition on the intersection. \footnote{Note that we are abusing notation by writing the interval $[0,1]$ as closed. These maps are colimits of \emph{collared} maps on open intervals $(0-\epsilon,1+\epsilon)$ as $\epsilon\to 0$.} Moreover, the last two maps must interpolate between the transition data for the endpoints. First observe that since $U$ is a convex open subset of $\RR^k$, for some $k$, $U$ smoothly deformation retracts to any fixed point $u\in U$. Choose such a retraction $r_t:U\to U$ which is constant on the collars of the interval. Precomposing all the cocycle maps determining $f$ with this deformation retraction gives an edge connecting $f$ to a map which is constant on the factor $U$. By composing edges, we can therefore assume that $f$ is constant on $U$.
 
Let $f^{\beta}:U\times W_{\beta}\to D^N/\partial D^N$ be the the cocycle maps determined by $d_0H=f$. Let $\{f^{\gamma}\}$ be the subset of all $\{f^{\beta}\}$ factoring through $D^N$. Since the chosen cover is countable, in fact finite, Lemma \ref{rvd} implies that we can find a small sphere in $D^N$ containing regular values $y_{\gamma}$ for each $f^{\gamma}$. 
For each $\gamma$ in this finite set, let $H^{\gamma}$ be post-composition with the smooth family of diffeomorphisms guaranteed by Lemma 
 \ref{lem sfd}. Choose $g^{\alpha\beta}$ and $g^{Z\beta}$ to be the transition maps of $f$, constant in the direction of the interval $[0,1]$. For all the $f^{\beta}$'s which do not factor through $D^N$, choose $H^{\beta}=f^{\beta}$ and $H^Z=f^{Z}$, constant in the direction of $[0,1]$. The edge corresponding to these choices has the desired property.
 
Finally, it is clear from the definition that if $f$ is constant on $U$, then $g=d_1H$ is constant on $U$ and for each $\beta$, the projection of the preimage $N_{\beta}=g^{-1}(0)\to U$ is a trivial bundle $N_{\beta}\cong M_{\beta}\times U$ with $M_{\beta}=(g\circ u)^{-1}(0)$ and $u:W_{\beta}\to U\times W_{\beta}$ induced by a point $u:\ast\to U$.
\endofproof

With all the above setting up, we are now ready for our first main theorem. 
\begin{proposition}
[Sheaves of connected components]
The map \eqref{abstract collapse map} is an isomorphism.
\end{proposition}
\theproof
We construct the inverse map to ${\bf PT}$ objectwise. Suppose we are given a map
\(\label{inv pi0}
\xymatrix{
f^{\prime}:D^{d+N-1}/\partial D^{d+N-1}\wedge U_+
\ar[r] & 
{\rm Th}\big(\mathcal{U}^{\perp}(d,N)\big)
}
\)
representing an element in the right hand side of \eqref{abstract collapse map}. We want to construct a convenient representative for the homotopy pullback
$$
\xymatrix@R=1.5em{
D^{d+N-1}/\partial D^{d+N-1}\wedge U_+ \ar[rrr]^-{f^{\prime}} &&& 
{\rm Th}\big(\mathcal{U}^{\perp}(d,N)\big)
\\
N \ar[rrr]\ar[u] &&& 
{\rm Gr}(d,N)\ar[u]_-{\rm zero\ section}\;.
}
$$
Since ${\rm Th}(\mathbfcal{U}^{\perp}(d,N))\simeq {\rm Th}(\mathcal{U}^{\perp}(d,N))$ and ${\bf Gr}(d,N)\simeq {\rm Gr}(d,N)$, it suffices to consider instead the homotopy pullback diagram
\(
\label{eq conv rep stack}
\xymatrix@R=1.5em{
D^{d+N-1}/\partial D^{d+N-1}\wedge U_+ \ar[rrr]^-{f^{\prime}} &&& 
{\rm Th}\big(\mathbfcal{U}^{\perp}(d,N)\big)
\\
N \ar[rrr]\ar[u] &&& 
{\bf Gr}(d,N)\ar[u]_-{\rm zero\ section}\;.
}
\)
Let $\{W_{\beta}\}\cup \{Z\}$ be a cover of $D^{d-1+N}/\partial D^{d-1+N}$, as in 
 Lemma \ref{Lem cover}. 
 We will define the inverse map locally on the covering and show that we can ``glue" the results 
 together to get the desired object. For the sake of readability, we will break the proof up into steps.
 
 \medskip
 \noindent {\bf Step 1} {\bf (The local calculation)} For each $W_{\beta}$, let $N_{\beta}$ be defined by the 
 homotopy pullback
$$
\xymatrix@R=1.5em{
W_{\beta}\times U\ar[rrr]^-{f^{\prime}_{\beta}} &&& 
{\rm Th}\big(\mathbfcal{U}^{\perp}(d,N)\big)
\\
N_{\beta} \ar[rrr]\ar[u] &&& 
{\bf Gr}(d,N)\ar[u]_-{\rm zero\ section}\;,
}
$$
where $f_{\beta}$ is the restriction to $W_{\beta}$. Since $W_{\beta}$ is zero-truncated, we 
can use the explicit presentations of the stacks on the right discussed in Section \ref{Sec gras}
 to compute this homotopy pullback explicitly as the strict pullback
$$
\xymatrix@=1.5em{
W_{\beta}\times U\ar[rr]^-{f^{\prime}_{\beta}} && D^{N}/\partial D^N
\\
N_{\beta} \ar[rr]\ar[u] &&\ast \ar[u]_-{0}
}\;,
$$
where the vertical map is induced by the map $0:\ast\to D^N\to D^N/\partial D^N$ which picks out 
zero in the closed unit disc. Since $W_{\beta}$ is representable, the definition of $D^N/\partial D^N$ implies that $f^{\prime}_{\beta}$ can be identified with an equivalence class smooth functions $f^{\prime}_{\beta}:W_{\beta}\times U\to D^N$. This equivalence identifies all smooth functions which factor through the boundary. Notice that if $f^{\prime}_{\beta}$ factors through the boundary, then since the boundary is far away from zero, $N_{\beta}=\emptyset$ and this is true for any representative of the class. Thus, we can assume that $f^{\prime}_{\beta}$ factors as $f^{\prime}_{\beta}:W_{\beta}\times U\to D^N\to D^N/\partial D^N$ and we are reduced to the pullback square of smooth manifolds
\(\label{rest f}
\xymatrix@=1.5em{
W_{\beta}\times U\ar[rr]^-{f^{\prime}_{\beta}} && D^{N}
\\
N_{\beta} \ar[rr]\ar[u] &&\ast \ar[u]_-{0}
\;.
}
\)
If $f^{\prime}_{\beta}$ has zero as a regular value, then $N_{\beta}$ will be a smooth submanifold of $W_{\beta}$. Here is where we will use the fact that the inverse need only be well-defined 
at the level of $\widetilde{\pi}_0$, and our transversality from 
Proposition \ref{homotopy edge} implies that we can assume this is the case. Moreover (by the same proposition), we can assume that $f_{\beta}$ is constant in the $U$ direction so that $N_{\beta}\cong M_{\beta}\times U$. 

\medskip
\noindent {\bf Step 2.} {\bf (Understanding the construction as a ``local" inverse)}
Each $N_{\beta}$ comes as the preimage of a regular value of $0\in D^{N}$. Hence the kernel of $df_{\beta}$ is identified with the tangent bundle $TN_{\beta}\cong TM_{\beta}\times \underline{\RR^{{\rm dim}(U)}}$, where the second factor is the trivial bundle of rank ${\rm dim}(U)$ over $U$.  In addition, $N_{\beta}$ inherits a local framing of the normal bundle from a choice of basis of orthonormal basis $\{v_i\}$, with $v_i\in D^N$. Summarizing in our case, the normal bundle is identified with $T^{\perp}N_{\beta}\cong T^{\perp}M_{\beta}\times U$ where the $T^{\perp}M_{\beta}$ is normal to $TM_{\beta}$ in $D^{d+N-1}$ and $f^{\prime}_{\beta}$ sends each framing of $T^{\perp}M_{\beta}\times U$ to a basis $\{v_i\}$. 

The map $f^{\prime}_{\beta}$ appearing in \eqref{rest f} is implicitly the restriction of the original map $f^{\prime}$ to the factor $D^N/\partial D^N$. Consequently, 
each $f^{\prime}_{\beta}$ also determined a map to ${\rm V}(N,d)$, which chooses a $d$-plane in $\RR^{d+N}$, up to some ambiguity parametrized by a choice of orthonormal basis of the complement. The discussion above shows that each framing of $T^{\perp}M_{\beta}$ gives a basis of the complement of this $d$-plane. After dividing out by the action of ${\rm O}(N)$, we thus see that $f^{\prime}_{\beta}$ gives exactly the map $\gamma$ in \eqref{maps to universal}, but restricted to some local patch of $M$. Thus, provided we can glue the $N_{\beta}\cong M_{\beta}\times U$ together, we will immediately see that this 
construction gives a two-sided inverse.

\medskip
\noindent {\bf Step 3.} {\bf (Gluing)} It remains to show that the manifolds $N_{\beta}$ glue properly. Recall that, since $\{W_{\beta}\times U\}\cup \{Z\wedge U_+\}$ 
is a cover of $D^{d+N-1}/\partial D^{d+N-1}\wedge U_+$, a map of the form \eqref{inv pi0} is uniquely 
determined by the data of a collection of maps 
$f^{\prime}_{\beta}:W_{\beta}\times U\to D^N/\partial D^N\times {\rm V}(N,d)$ on open sets and maps 
$g_{\beta\gamma}:W_{\beta\gamma} \times U\to {\rm O}(N)\times {\rm V}(N,d)$ on intersections satisfying 
some cocycle conditions. We will only be concerned with the condition that affects the base level. That is,
 we have $f^{\prime}_{\beta}=g_{\beta\gamma}f^{\prime}_{\gamma}$
on intersections, where the juxtaposition on the right is given by the ${\rm O}(N)$-action. Since 
$0\in D^{N}/\partial D^N$ is fixed by the ${\rm O}(N)$-action, we have
$$
(W_{\beta}\times U)\cap N_{\gamma} = f_{\gamma}\vert_{W_{\beta\gamma}}^{-1}(0) =  
(g_{\beta\gamma}f^{\prime}_{\beta})\vert_{W_{\beta\gamma}}^{-1}(0) \;.
$$
Since the $g_{\beta\gamma}$ act by invertible linear maps, we immediately get 
$$
(g_{\beta\gamma}f^{\prime}_{\beta})\vert_{W_{\beta\gamma}}^{-1}(0) =(f^{\prime}_{\beta})^{-1}\vert_{W_{\beta\gamma}}(0)=(W_{\gamma}\times U)\cap N_{\beta}\;.
$$
 Then the union $\bigcup_{\beta}N_{\beta}=:N$ defines a smooth submanifold of 
 $(D^{d+N-1})^{\rm o}\times U\subset D^{d+N-1}/\partial D^{d+N-1}\wedge U_+$. The projection $N\to U$ is a fiber bundle with typical fiber $M=\bigcup_{\beta} M_{\beta}$ and we
have determined a submanifold $N\subset D^{d+N-1}\times U \subset \RR^{d+N-1}$, which is diffeomorphic to the product $M\times U$ and this defines an object in $\BB\mathscr{C}{\rm ob}_d(U)$. By the local explanation above (Step 2), we immediately see that this map is a two-sided inverse.
\endofproof

Finally, we observe that since  both $\BB\mathscr{C}{\rm ob}_d$ and $\BB{\rm Conc}^{>}\big(\Omega^{\infty-1}_{{}_{D^1/\partial D^1}}{\MM\TT}(d)\big)$ are zero-truncated, they are equivalent to their sheaves of connected components. This immediately gives the first of our main theorems 
\begin{theorem}\label{theorem 1}
The map \eqref{abstract collapse map} induces an equivalence
$$
\xymatrix{
{\bf PT}:\BB\mathscr{C}{\rm ob}_d\ar[r]^-{\simeq} & \BB{\rm Conc}^{>}\big(\Omega^{\infty-1}_{{}_{D^1/\partial D^1}}{\MM\TT}(d)\big)
\;.}
$$
\end{theorem}

\section{Higher smooth tangential structures} 
\label{Sec higher}

\subsection{Adding $G$-structure}
\label{Sec G}

In this section, we describe how to add smooth tangential structures to the picture. 
Much of the machinery established in the previous section can be used to include these various structures 
on the tangent bundle via pullback. In general, if we are given a Lie group $G$ and a faithful  representation $G\into {\rm O}(d)$, then delooping this map gives rise to a map  
$$
\theta:\BB G\longrightarrow \BB{\rm O}(d)\;.
$$
For example, $\theta$ could be induced by a representation of $G$.  Note that although it is not true that ${\bf Gr}(d,N)\to \BB{\rm O}(d)$ as $N\to \infty$, there is still a map ${\bf Gr}(d,\infty)\to \BB{\rm O}(d)$ which is induced by taking the model ${\bf Gr}(d,\infty)={\rm V}(d,\infty)/\!/{\rm O}(d)$ and projecting out ${\rm V}(d,\infty)$. These observations lead to the following definition, in analogy to the tangential structures considered in \cite{GMTW}.
\begin{definition}
[Smooth Grassmannian stack and universal bundle with $\theta$-structure]
{\bf (i)} We define the {\rm smooth Grassmanian stack with $\theta$-structure} as the homotopy pullback
$$
\xymatrix@R=1.5em{
{\bf Gr}(\theta,N)\ar[rr]\ar[d]_{\theta_{N}}&& \BB G\ar[d]^{\theta}
\\
{\bf Gr}(d,N)~\ar@{^{(}->}[r] & {\bf Gr}(d,\infty)\ar[r] & \BB {\rm O}(d)\;.
}
$$
\noindent {\bf (ii)} We define the {\rm universal bundle}
 $\mathbfcal{U}(\theta,N)\to {\rm Gr}(\theta,N)$ as the homotopy pullback 
$$
\xymatrix@=1.5em{
\mathbfcal{U}(\theta ,N)\ar[rr]\ar[d]&& \mathbfcal{U}(d,N)\ar[d]
\\
{\bf Gr}(\theta,N)~\ar@{^{(}->}[rr]^-{\theta_{N}} && {\bf Gr}(d,N)\;.
}
$$
\noindent {\bf (iii)} We define the {\rm universal orthogonal complement bundle}
 $\mathbfcal{U}^{\perp}(\theta,N)\to {\bf Gr}(\theta,N)$ as the homotopy pullback 
$$
\xymatrix@=1.5em{
\mathbfcal{U}^{\perp}(\theta ,N)\ar[rr]\ar[d]&& \mathbfcal{U}^{\perp}(d,N)\ar[d]
\\
{\bf Gr}(\theta,N)~\ar@{^{(}->}[rr]^-{\theta_{N}} && {\bf Gr}(d,N)\;.
}
$$
\end{definition}

Since these stacks are given by pullbacks, they enjoy much of the same properties as their 
counterparts without $\theta$-structure. In particular, we have an extension of the splitting result 
 in Proposition \ref{pullback trivializes} to include $\theta$-structures. 
 
\begin{proposition}\label{decomposition}
The pullback of $\mathbfcal{U}^{\perp}(\theta,N+1)$ to ${\bf Gr}(\theta,N)$ decomposes as 
the sum $\mathbfcal{U}^{\perp}(\theta,N)\oplus {\bf 1}$.
\end{proposition}
  \theproof
 Let $i:{\bf Gr}(\theta,N)\into {\bf Gr}(\theta,N+1)$ be the map induced by the inclusion on Grassmannian stacks. By Proposition \ref{pullback trivializes}, we have an isomorphism of bundles
\bea 
i^*\theta^*_{N+1}\mathbfcal{U}^{\perp}(d,N+1) &
\cong& \theta^*_{N+1}i^*\mathbfcal{U}^{\perp}(d,N+1)
\\
&\cong& \theta^*_{N+1}\big(\mathbfcal{U}^{\perp}(d,N)\oplus {\bf 1}\big)
\\
&\cong& \big(\theta^*_N\mathbfcal{U}^{\perp}(d,N)\big)\oplus {\bf 1}\;.
\eea

\vspace{-10mm}
 \endofproof
 
As in the case without $\theta$-structure, Proposition \ref{decomposition} 
gives rise to maps
 \(
 \label{DThTh}
 \xymatrix{
 D^1/\partial D^1\wedge {\rm Th}\big(\mathbfcal{U}^{\perp}(\theta,N)\big)
 \ar[r] &
  {\rm Th}\big(\mathbfcal{U}^{\perp}(\theta,N+1)\big)
  }\;,
 \)
 which turn ${\rm Th}(\mathbfcal{U}^{\perp}(\theta,N))$ into a $D^1/\partial D^1$-spectrum 
 as $N$ varies. We denote this spectrum accordingly as $\MM\TT(\theta)$. The homotopy pullbacks in Proposition \ref{decomposition} are easily seen to be zero-truncated. Indeed, fix a basepoint, say $x:\ast\to {\bf Gr}(\theta,N)$. Then we have a corresponding homotopy fiber sequence
 $$
 {\bf Gr}(\theta,N)\longrightarrow 
 \BB G\times {\bf Gr}(d,N)\longrightarrow 
 \BB{\rm O}(d)\;.
 $$
 Since $G\cong \widetilde{\pi}_1(\BB G)\into {\rm O}(d)\cong \widetilde{\pi}_1(\BB{\rm O}(d))$, the long exact sequence on sheaves of verifies the claim. It follows that each of the above homotopy pullbacks must be equivalent to its sheaf of components, which we denote by ${\rm Gr}(\theta,N)$, $\mathcal{U}(\theta,N)$ and $\mathcal{U}^{\perp}(\theta,N)$. The latter two objects are still vector bundles over ${\rm Gr}(\theta,N)$. We can,
  therefore, define the sheaf of bundle maps ${\rm Bun}(TM,\mathcal{U}(\theta))$ as the sheaf which on each Cartesian space $U$ assigns the set of bundle maps
 $$
 {\rm Bun}(TM,\mathcal{U}(\theta))(U):=\{p:TM\times U\to \mathcal{U}(\theta):p\ {\rm is\ a\ fiberwise\ bundle\ map}\}\;.
 $$
We can make the same definition using the stack $\mathbfcal{U}(\theta)$, and given that the canonical equivalence $q:\mathbfcal{U}(\theta)\overset{\simeq}{\to} \mathcal{U}(\theta)$ defines a bundle map, we have an equivalence ${\rm Bun}(TM,\mathbfcal{U}(\theta))\simeq {\rm Bun}(TM,\mathcal{U}(\theta))$.

\medskip
We define the smooth sheaf ${\rm Emb}_{\theta}(M,\RR^{d-1+\infty})$ as the pullback of sheaves
$$
\xymatrix@=1.5em{
{\rm Emb}_{\theta}(M,\RR^{d-1+\infty})\ar[rr]\ar[d] &&
 {\rm Bun}(TM,\mathcal{U}(\theta))\ar[d]
\\
{\rm Emb}(M,\RR^{d-1+\infty})\ar[rr]^-{TM} &&
 {\rm Bun}(TM,\mathcal{U}(d))\;,
}
$$
 where $TM$ denotes the canonical bundle map lifting the Gauss map $TM:M\times \RR \to {\rm Gr}(d,\infty)$. 
 
 \begin{remark}
 [Space vs. smooth sheaves of bundle maps]
 It is well known that the \emph{space} of bundle maps ${\rm Bun}(TM,\mathcal{U}(d))$ is contractible (see for example \cite[Lemma 5.1]{GMTW}). However, this is of course not the case as in smooth sheaves. Nevertheless, the geometric realization of ${\rm Emb}_{\theta}(M,\RR^{d-1+\infty})$ agrees with the corresponding space defined in \cite{GMTW}.
 \end{remark}
 
\medskip
In the same spirit as the previous section, we can form the smooth cobordism category 
with smooth $G$-structure (hence using classifying stacks ${\bf B}G$ rather than classifying 
spaces   ${\rm B}G$). We make the following definition.

 \begin{definition}  [Smooth cobordism category with $G$-structure]
  The {\rm cobordism category with $G$-structure} has as sheaf of objects and morphisms, respectively,  
 \bea
 {\rm Ob}(\mathscr{C}{\rm ob}_d)&:= &
 \coprod_{[M]}\RR\times {\rm Emb}_{\theta}(M,\; \RR^{d-1+\infty})/{\rm Diff}(M)\;,
 \\
  {\rm Mor}(\mathscr{C}{\rm ob}_d)&:= &
  \coprod_{[W]}\RR^{2}_+\times {\rm Emb}_{\theta}(W, \; [0,1]\times \RR^{d-1+\infty})/{\rm Diff}(W)\;.
\eea

\end{definition}
The sections of the sheaf of objects can be identified with a triple $(N,t,l)$, where $N\cong M\times U\subset \RR^{d+\infty-1}$ is a bundle of $(d-1)$-dimensional manifolds $M$, $t:U\to\RR$ is a smooth function and $l:N\cong M\times U\to \BB G$ is a lift of $TM:N\cong M\times U \to {\rm Gr}(d,\infty)\simeq {\bf Gr}(d,\infty)\to \BB{\rm O}(d)$. \footnote{Here we really mean a lift of any map in the connected component that $TM:N\to{\rm Gr}(d,\infty)$ corresponds to in $\map(N,{\bf Gr}(d,\infty))$.} The morphisms are identified similarly, with lifts $l_{t,t^{\prime}}:Z\cong W\times U\to \BB G$ which are required to restrict to the maps $l_t:N\to \BB G$ and $l_{t^{\prime}}:N^{\prime}\to \BB G$ on collars. 

\medskip
Notice that the proof of the theorem with $G$-structure is almost a direct extension of the proof without $G$-structure (Theorem \ref{theorem 1}), given by replacing $\mathbfcal{U}^{\perp}(d)$ with $\mathbfcal{U}^{\perp}(\theta)$. Indeed, we still have a well-defined collapse map
\(\label{PT map 2}
\xymatrix{
{\bf PT}:\widetilde{\pi}_0\big(\BB\mathscr{C}{\rm ob}_{\theta}\big)
\; \ar[r] & \;
\widetilde{\pi}_0\big(\BB{\rm Conc}^{>} \big(\Omega^{\infty-1}_{{}_{D^1/\partial D^1}}{\MM\TT}(\theta)
\big)\big)
}\;,
\)
defined via the universal property of the pullback, and the same proof of Proposition \ref{Prop SHG} applies. We will simply sketch how to construct the inverse map in this case. 
\begin{theorem}[Stacky Pontrjagin-Thom equivalence with $\theta$-structure]
\label{theorem 2}
The map \eqref{PT map 2} induces an equivalence
\(
\xymatrix{
{\bf PT}:\BB\mathscr{C}{\rm ob}_{\theta}
\; \ar[r] & \BB{\rm Conc}^{>} \big(\Omega^{\infty-1}_{{}_{D^1/\partial D^1}}{\MM\TT}(\theta)
\big)\big)
}
\)
\end{theorem}
\theproof
As discussed above, we will sketch the construction of the inverse to the map \eqref{PT map 2}. Fix a map
$$
\xymatrix{
f^{\prime}:D^{d+N-1}/\partial D^{d+N-1}\wedge U_+
\; \ar[r] & \;
{\rm Th}(\mathbfcal{U}^{\perp}(\theta,N))
}\;,
$$
and cover $D^{d+N-1}/\partial D^{d+N-1}$ as in the constructions in Section \ref{Sec PT}. 
We are again led to the homotopy pullback diagram 
\(\label{first square}
\xymatrix@=1.5em{
D^{d+N-1}/\partial D^{d+N-1}\wedge U_+
\ar[rr] && 
{\rm Th}(\mathbfcal{U}^{\perp}(\theta,N))
\\
M \ar[u]
\ar[rr] && 
{\bf Gr}(\theta,N)\ar[u]_-{\rm zero\ section}\;,
}
\)
where $M$ has yet to be identified. By definition, the stacks ${\rm Th}(\mathbfcal{U}^{\perp}(\theta,N))$ and ${\bf Gr}(\theta,N)$ fit into a pullback diagram
$$
\xymatrix@=1.5em{
 {\rm Th}(\mathbfcal{U}^{\perp}(\theta,N))\ar[d]\ar[rr] &&
   {\rm Th}(\mathbfcal{U}^{\perp}(d ,N))\ar[d]
\\
 {\bf Gr}(\theta,N)\ar[rr] &&  {\bf Gr}(d,N)
\;.
}
$$
Precomposing both vertical maps with the maps induced from the zero section, we get a diagram
$$
\xymatrix@=1.5em{
{\bf Gr}(\theta,N)\ar[rrr]\ar[d] &&&  {\rm Th}(\mathbfcal{U}^{\perp}(\theta,N))
\ar[d]\ar[rr] && \ar[d] {\bf Gr}(\theta,N)
\\
 {\bf Gr}(d,N)\ar[rrr]^-{\rm zero\ section} &&& 
 {\rm Th}(\mathbfcal{U}^{\perp}(d ,N)) \ar[rr] &&  {\bf Gr}(d,N)
\;.
}
$$
The outer diagram has top and bottom horizontal maps as identities and is, therefore, homotopy Cartesian. The right diagram is homotopy Cartesian by definition and so the Pasting Lemma implies that the left square is homotopy Cartesian. Combining the left square above with diagram \eqref{first square}, we see that the pasting law for homotopy pullbacks allows us to compute $M$ as the homotopy pullback of the resulting outer square, which was already computed in Theorem \ref{theorem 1}. This presents $M$ manifestly as a smooth manifold, equipped with all the necessary data to identify it with an object in ${\rm Cob}_{\theta}$.
\endofproof

This then gives a Pontrjagin-Thom construction at the level of $G$-bundles with 
data of gauge transformations. The group $G$ can be taken to be a subgroup of ${\rm O}(n)$,
such as ${\rm U}(m)$, ${\rm Sp}(m)$, or even exceptional groups if $n$ is large enough.

\subsection{Adding geometric data}
\label{Sec geom}

We now explain how to 
add geometric data, such as metrics and connections. 
Let us fix a sheaf ${\cal F}$ which is equipped with an ${\rm O}(d)$-action. Then we can add ${\cal F}$-structure by considering the relevant pullbacks along the forgetful maps ${\rm O}(d)/\!/{\cal F}=:\BB{\rm O}(d)_{\cal F}\to \BB{\rm O}(d)$, induced by the projection ${\cal F}\to \ast$. We name the resulting smooth stacks ${\bf Gr}(d,N)_{\cal F}$, $\mathbfcal{U}(d,N)_{\cal F}$ and $\mathbfcal{U}^{\perp}(d,N)_{\cal F}$ accordingly. Again each smooth stack is zero-truncated and equivalent to a smooth sheaf ${\rm Gr}(d,N)_{\cal F}$, $\mathcal{U}(d,N)_{\cal F}$, and $\mathcal{U}^{\perp}(d,N)_{\cal F}$ and all the properties that hold for the stacks without ${\cal F}$-structure 
still hold for their ${\cal F}$ counterparts, In particular, we have maps (cf. maps 
\eqref{DThTh})
$$
\xymatrix{
D^1/\partial D^1\wedge {\rm Th}(\mathbfcal{U}^{\perp}(d,N)_{\cal F}) 
\; \ar[r] & \; {\rm Th}\big(\mathbfcal{U}^{\perp}(d,N+1)_{\cal F}\big)
}\;,
$$
and therefore a $D^1/\partial D^1$-spectrum $\MM \TT(d)_{\cal F}$. Similarly, we have maps 
$$
\xymatrix{
D^1/\partial D^1\wedge {\rm Th}(\mathcal{U}^{\perp}(d,N)_{\cal F}) 
\; \ar[r] & \; {\rm Th}\big(\mathcal{U}^{\perp}(d,N+1)_{\cal F}\big)
}\;,
$$
which gives an equivalent $D^1/\partial D^1$-spectrum $MT(d)_{\cal F}$. On the bordism side, 
we can define the sheaf of bundle maps ${\rm Bun}(M,\mathcal{U}(d)_{\cal F})$ which lift $TM$ as in the case of $G$-structure. We define
$$
\xymatrix@=1.5em{
{\rm Emb}^{\cal F}(M,\RR^{d-1+\infty})\ar[rr]\ar[d] && 
{\rm Bun}(TM,\mathcal{U}(d)_{\cal F})\ar[d]
\\
{\rm Emb}(M,\RR^{d-1+\infty})\ar[rr]^-{TM} &&
 {\rm Bun}(TM,\mathcal{U}(d))\;,
}
$$
and define the smooth cobordism category ${\rm Cob}_{d}^{\cal F}$ accordingly. 

\begin{definition}
[Smooth cobordism category with geometric data]
Define the \emph{smooth cobordism category with geometric data 
$\mathscr{C}{\rm ob}_{d}^{{\cal F}}$}
 by
 \bea
 {\rm Ob}(\mathscr{C}{\rm ob}_{d}^{{\cal F}})&:=&
 \coprod_{[M]}{\rm Emb}^{\cal F}(M,\RR^{d-1+\infty})/{\rm Diff}(M)\;,
 \\
 {\rm Mor}(\mathscr{C}{\rm ob}_{d}^{{\cal F}})&:=&\coprod_{[W]}{\rm Emb}^{\cal F}(W,\RR^{d-1+\infty})/{\rm Diff}(W)\;.
 \eea
 \end{definition}
 In this case, the sections of the sheaf of objects of the resulting smooth cobordism category are triples  $(N,t,l)$, with $N\to U$ a bundle of $(d-1)$-dimensional manifolds, a smooth function $t\in C^{\infty}(U;\RR)$ and a lift $l:N\cong M\times U\to \BB {\rm O}(d)_{\cal F}$ of $TM$, which is natural in $U$. The bordisms are identified similarly with the usual compatibility on collars. The definition and construction are very versatile and allow for geometric data such as the following. 
 
 \begin{example} [Connections]
 Let ${\cal F}=\nabla:=\Omega^1(-;\mathfrak{o}(d))$ and let ${\rm O}(d)$ act on $\nabla$ via gauge transformations. Then the stack $\BB {\rm O}(d)_{\nabla}$ is the moduli stack of connections on $M$. 
 The corresponding bordism category is the category with objects manifolds equipped with connections 
 on the tangent bundle. The morphisms are bordisms equipped with connections on the tangent bundle 
 which extend those of the bounding manifolds up to gauge transformation.
 \end{example} 
\begin{example} [Metrics] Let ${\cal F}=\mathcal{M}(\RR^d):=C^{\infty}(-,{\rm Sym}_d\cap GL_d)$ 
be the sheaf of metrics on $\RR^d$. For each Cartesian space $U$, $\mathcal{M}(\RR^d)(U)=\{ {\rm Metrics\ on\ the\ bundle}\ p:\RR^d\times U\to U\}$. Then ${\rm O}(d)$ acts on $\mathcal{M}(\RR^d)$ via 
conjugation of symmetric matrices. The stack $\BB {\rm O}(d)_{{\cal M}(\RR^d)}$ is the moduli stack 
of metrics on $M$. The corresponding bordism category has objects -- smooth manifolds equipped with 
Riemannian metrics, and morphisms -- bordisms equipped with metrics which extend those of the 
bounding manifolds, up to local change of basis.
\end{example}

Just as in the case of $G$-structure, the collapse map is defined using information in the modified cobordism category and one observes that the introduction of the sheaf ${\cal F}$ causes no difficulty. In fact, the same proof of Theorem \ref{theorem 2} works in this case. This leads us to the following theorem. 
\begin{theorem}
[The Pontrjagin-Thom equivalence with geometric structure]
\label{final thm}
We have an equivalence of smooth stacks
$$
\xymatrix{
{\bf PT}:\BB{\rm Cob}_{d}^{{\cal F}}
\; \ar[r]^-{\simeq} & \;
 \BB{\rm Conc}^>(\Omega^{\infty-1}_{D^1/\partial D^1}{\MM\TT}(d)_{\cal F})
}\;.
$$
\end{theorem}

We now quickly produce a nontrivial example already in the case of cobordism of 
zero-dimensional manifolds.

\begin{example}
[An application for $d=1$]
 Let ${\cal F}:=\Omega^1_{\rm cl}$ be the sheaf of closed 1-forms on Cartesian spaces. The ${\rm O}(1)=\ZZ/2$ acts on this sheaf via gauge transformations, which in this case is trivial. We readily calculate
$$
\mathcal{U}^{\perp}(1,N)_{\cal F}\cong \Omega^1_{\rm cl}\times \mathcal{U}^{\perp}(1,N)\;,
$$
so the corresponding Thom stack simplifies to ${\rm Th}(\mathcal{U}^{\perp}(1,N)_{\cal F})\cong \big(\Omega^1_{\rm cl}\big)_+\wedge {\rm Th}(\mathcal{U}^{\perp}(1,N))$. The geometric realization of the smooth sheaf $\Omega^1_{\rm cl}$ can be identified with $K(\RR,1)$ \cite{Pa}. Then, we see that if we geometrically realize the smooth stack $\BB{\rm Conc}^>(\Omega^{\infty-1}{\MM\TT}(0)_{\cal F})$, we get the space
\bea
\vert \BB{\rm Conc}^>(\Omega^{\infty-1}{\MM\TT}(0)_{\cal F})\vert &\simeq&
 \vert \Omega^{\infty-1}{\MM\TT}(0)_{\cal F} \vert 
 \\
& \simeq & \Omega^{\infty-1}\big(K(\RR,1)_+\wedge (\RR P^{\infty})^{-L}\big)\;,
\eea
where $L\to \RR P^{\infty}$ is the canonical line bundle with virtual inverse $-L$, and $(\RR P^{\infty})^{-L}$ denotes the corresponding Thom spectrum. Notice, however, that this is the suspension of the sphere spectrum  $(\RR P^{\infty})^{-L}\simeq \Sigma\mathbb{S}$. Indeed, the unit disc bundle $D\to \RR P^{n+1}$ can be obtained by intersecting each hyperplane with the unit disc $D^{n+1}\subset \RR^{n+1}$. So the Thom space is readily identified with $D^{n+1}/\partial D^{n+1}\simeq S^{n+1}$.

On the bordism side, we see that the sheaf of objects, evaluated at a cartesian space $U$, are zero-dimensional smooth manifolds 
$M=\{x_i\}_{i=1}^n\times U \subset \RR^{\infty-1}\times U$, equipped with maps $M\times U\times \RR \to \Omega^1_{\rm cl}\times \BB\ZZ/2$, lifting 
$$
\xymatrix{
T^\pi M:M\times U\times \RR \; \ar[r] &  \;
{\rm Gr}(1,\infty)\simeq \RR P^{\infty}\simeq S^{\infty}/\!/\ZZ/2
\; \ar[r] & \;
 \ast/\!/\ZZ/2
 }\;.
$$
Such lifts are in bijective correspondence with closed 1-forms, thought of as flat connections on the the line bundle $M\times U\times \RR\to M\times U$. The morphisms are embedded 1-dimensional manifolds $W\times U\subset \RR^{\infty}\times U$, where $W$ is a collection of $n$ smooth paths with boundary the disjoint union $M_0\sqcup M_1$, $M_0\cong M_1$. Graphically,
we have
$$
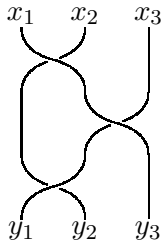
\vcenter{\xy 0;/r2pc/:
     (0,.2)*{x_1}, (1,.2)*{x_2}, (2,.2)*{x_3}, (0,-3.2)*{y_1}, (1,-3.2)*{y_2}, (2,-3.2)*{y_3}, 
    @={(0.5,-0.5),(0,0)}, @@{*{\vtwistneg}},(0,-2),\vtwist,
    @i@={(0,-1),(2,0),(2,-2)} @@{="save";"save"-(0,1),**@{-}}
    \endxy}
$$
We also have to consider maps $W\times U\to \Omega^1_{\rm cl}\times \BB\ZZ/2$, which restrict to the given forms on collared neighborhoods of the boundary. A closed 1-form on $W$ is equivalent to associating a nonzero real number on each strand, multiplication by which gives the notion of parallel transport along the strand. The group generated under composition is then seen to be the group generated by elementary matrices: $GL_n(\RR)$. Since this group is discrete, the lifts which restrict fiberwise to lifts of $T^\pi M$ are constant in the direction of $U$. Therefore, after geometrically realizing, we get
$$
\vert \BB\mathscr{C}{\rm ob}_{d}^{{\cal F}}
\vert\simeq {\rm B}\Big(\coprod_{n\geq 0}{\rm BGL}_n(\RR)\Big)\;.
$$
After looping, we get the group completion (see \cite{Segal} \cite{Ad2})
and this gives the identification 

$$
\ZZ\times {\rm BGL}_{\infty}(\RR)^+
\simeq \Omega^{\infty}\Sigma^{\infty}K(\RR,1)_+\;.
$$
This technique of identifying certain infinite loops spaces via geometric realization of our equivalence deserves a fuller development and indeed we will revisit it elsewhere. 
\end{example}

\medskip
Lastly, we can combine Theorem \ref{theorem 2} ($\theta$-structures) 
and Theorem \ref{final thm} (${\cal F}$-structures)
to give Theorem \ref{Thm most} in the Introduction.

\end{document}